\documentclass[11pt,notitlepage,twoside]{article}
\usepackage{latexsym}
\usepackage{amsmath}
\usepackage{amsfonts}
\usepackage{amssymb}

\usepackage{xcolor,rotating,epic,eepic}\usepackage{color}

\newcommand{\e }{\varepsilon }

\renewcommand{\i }{\iota}

\newcommand{\n }{\nabla }

\newcommand{\be}{\begin{equation}}
\newcommand{\ee}{\end{equation}}
\newcommand{\ba}{\begin{align*}}
\newcommand{\ea}{\end{align*}}
\newenvironment{pf}{\noindent{\bf Proof.}\enspace}{
\hfill$\Box$\medskip}
\newenvironment{pfn}[1]{\noindent{\bf Proof of {#1}\enspace}}{
\hfill$\Box$\medskip}
\newcommand{\R}{\mathbb{R}}
 
\newcommand{\Z}{\mathbb{Z}}

\newcommand{\N}{\mathbb{N}}

\newtheorem{theorem}{Theorem}[section]
\newtheorem{definition}{Definition}[section]
\newtheorem{lemma}{Lemma}[section]
\newtheorem{proposition}{Proposition}[section]

\newtheorem{corollary}{Corollary}[section]

\let \n = \noindent

\author{  Zakaria Boucheche \\
{\footnotesize \noindent  Laboratory of applied mathematics and harmonic analysis LR11ES52, Sfax, Tunisie}}

\title { \Large \textbf{Existence result under flatness condition for a nonlinear elliptic
 equation with Sobolev exponent}}

\begin{document}

\date{ }

\maketitle {\footnotesize} \noindent {\bf Abstract.} In this paper,
we consider the following nonlinear elliptic equation with Dirichlet
boundary condition: $-\Delta u=K(x)u^{\frac{n+2}{n-2}},\,  u>0$ in
$\Omega,\, u=0$ on $\partial\Omega$, where $\Omega$ is a smooth
bounded domain in $\mathbb{R}^n,$ $n\geqslant 4,$ and $K$ is a
$\mathcal{C}^1$-positive function in $\bar{\Omega}$. Under the
assumption that the order of flatness at each critical point of $K$
is $\beta \in ]\,n-2,\,n[,$ we give precise estimates on the looses
of the compactness, and we prove an existence result through an
Euler-Hopf type formula.\\\\\noindent{ \footnotesize { {\bf Key
words :}\,\,Elliptic equation, critical Sobolev exponent, loss of
compactness, variational method, critical points at infinity}}\\
\noindent{ \footnotesize { {\bf Mathematic subject
classification:}\,\,35J60}}
\section{Introduction and main result }
\def\theequation{1.\arabic{equation}}\makeatother
\setcounter{equation}{0} In this work, we look for solution for the
following nonlinear problem  under the Dirichlet boundary condition
\begin{equation}\label{problem} \left\lbrace
\begin{aligned}
-\Delta u&=K(x)u^{\frac{n+2}{n-2}} \quad \qquad \quad \mbox{in}\,\Omega\\
u&>0 \,\,\,\,\,\,\,\,\,\,\,\,\qquad \quad \qquad \quad \mbox{in}
\,\Omega\\
u&=0 \,\,\,\,\,\,\,\,\,\,\,\qquad \quad \, \qquad \quad \mbox{on}
\,\partial\Omega,
\end{aligned}
\right.
\end{equation}
where $\Omega$ is a bounded smooth domain of
 $\mathbb{R}^n$, $n\geqslant 4$, and $K$ is a
 $\mathcal{C}^1$-positive function in
$\bar{\Omega}.$\\

 One motivation to study this equation comes from its resemblance
to the prescribed scalar curvature problem in conformal geometry,
which consists on finding suitable conditions on a given function
$K$ defined on $M$ to be the scalar curvature of a metric
$\widetilde{g}$ conformally equivalent to $g$, where $(M,g)$ is an
$n$-dimensional Riemannian manifold without boundary. The special
nature of the problem (\ref{problem}) appears when we consider it
from the variational viewpoint. Indeed, although this problem enjoys
a variational structure in the sense that its solutions can be
interpreted as critical points of some functional, its associated
Euler-Lagrange functional does not satisfy the Palais-Smale
condition. This means that there exist noncompact sequences along
which the functional is bounded and its gradient goes to zero. This
is due to the non compactness of the embedding $H^1_0(\Omega)$ into
$L^{\frac{2n}{n-2}}(\Omega).$\vskip 2 mm \hrulefill{} \vskip 2 mm
{\fontfamily{ppl}\fontseries{b}%
E-mail address: Zakaria.Boucheche@ipeim.rnu.tn} \vskip 4 mm In the
case of manifolds without boundary, this problem has been widely
studied in various works, see, for example, the monograph \cite{A}
and the references therein. In contrast to the extensive literature
regarding the prescribed scalar curvature problem on manifolds
without boundary, in particular on spheres, there are few known
results of (\ref{problem}); see, for example, \cite{H2} and
\cite{14} for $n=4,$ \cite{Hebey} for $n> 4,$
$\Omega$ is a ball and $K=K(\|x\|),$ \cite{14} for $n\geq 4.$\\

One group of existence results have been obtained under hypotheses
involving $\Delta K$ at the critical points $y$ of $K.$ For example,
in \cite{14} it is assumed that $K$ is a Morse function
and\begin{equation}\label{nd}\Delta K(y)\neq 0,\,\,\,\forall\,\,y\in
\mathcal{K},\,\,\forall\,\,y\in \{\,x\in \Omega\,/\,\nabla
K(x)=0\,\}.\end{equation}Under the condition (\ref{nd}), an
Euler-Hopf criterion for $K$ was provided to find solution for the
problem (\ref{problem}). In \cite{14} we were explained that this
global criterion is not satisfied for higher dimensional case $n\geq
5.$ Naturally one may ask a similar question for all dimensions
$n\geq 4,$ namely, which function $K(x)$ on $\bar{\Omega}$ arise a
solution for the problem (\ref{problem}) under a global Euler-Hopf
criterion?\\

In order to state our main result, we need to introduce some
notations, and state the assumptions that we are using in our paper.
We denote by $G$ the Green's function and by $H$ its regular part,
that is for each $x\in\Omega$,
\begin{equation}\label{c0}
\left\lbrace
\begin{aligned}
G(x,y)&=|x-y|^{-(n-2)}-H(x,y)\quad \qquad \quad \mbox{in}\, \Omega\\
\Delta H(x,\,.\,)&=0 \qquad \qquad \qquad \quad \,\,\,\qquad \quad \qquad \quad \mbox{in} \,\Omega\\
G(x,\,.\,)&=0\,\,\,\qquad \qquad \qquad \qquad \qquad \qquad \quad
\mbox{on}\,\partial\Omega.
\end{aligned}
\right.
\end{equation} Let $K:\,\bar{\Omega}\rightarrow \R$ be a $\mathcal{C}^1$ positive
function.\\
$\mathbf{(A_1)}$ Assume that, for each $x\in\partial\Omega$, we have
$  \displaystyle\frac{\partial K(x)}{\partial \nu}<0 , $ where $\nu$
is the outward normal vector on $\partial\Omega$.\\Let
$$\mathcal{K}:=\{\,x\in \Omega\,/\,\nabla K(x)=0\,\}$$ the set of the
critical points of $K$ in $\Omega.$\\\\Throughout this paper, we
assume that $K$ satisfies the following flatness
condition:\\\\
$\mathbf{(f)}_\beta$ for each critical point $y$ of $K,$ there exist
$\beta:=\beta(y) \in ]n-2,\,n[,\,\,\mathrm{and\,\,}\eta > 0$ such
that in some local coordinates system centered at $y,$ we
have$$K(x)=K(y)+\sum_{i=1}^n\,b_k\bigl|(x-y)_k\bigr|^\beta\,+R(x),\,\,\forall\,\,x\in
B(y,\,\eta),$$where $b_k:=b_k(y)\neq
0,\,\,\forall\,\,k=1,\,\dots,\,n,$ and $R(z)$ is $\mathcal{C}^1$
near $0$ with$$\lim_{z\rightarrow
0}|R(z)||z|^{-\beta}=0\,\,\mathrm{and\,\,}\lim_{z\rightarrow
0}|\nabla R(z)||z|^{1-\beta}=0.$$Notice the
following.\begin{description}
\item[$\bullet$] The $\mathbf{(f)}_\beta$-assumption
was used widely as a standard assumption to guarantee the existence
of solution to the scalar curvature problem on closed manifolds.
However, for technical reason, it is assumed more regularity for the
function $R$ near $0,$ with the following
assumption$$\lim_{z\rightarrow 0}\sum_{s=0}^{[\beta]}|\nabla^s
R(z)||z|^{s-\beta}=0,$$where $[\beta]$ denotes the integer part of
$\beta;$ see, for example, \cite{Li}. Thus, the
$\mathbf{(f)}_\beta$-assumption mentioned above can be seen as a
refined condition.
\end{description}For each $y_i\in \mathcal{K},$
we will denote, if necessary, by $\beta_i$ for its order of
flatness.\\
For each $s-$tuple, $s\geq 1,$ of distinct points
$\tau_s:=(y_{i_1},\,\dots,\,y_{i_s})$ such that $y_{i_k}\in
\mathcal{K},\,\,\forall\,\,k=1,\dots,s,$ we define a $s\times s$
symmetric matrix $M(\tau_{s})=(m_{ij})$ by
\begin{equation}\label{matrice}
\begin{aligned}&m_{jj}:=\frac{H(y_{i_j},\,y_{i_j})}{\bigl(K(y_{i_j})\bigr)^{\frac{n-2}{2}}},\,\,
\forall\,\,1\leq j\leq s
\\&{}m_{jk}:=-\frac{G(y_{i_j},\,y_{i_k})}{\bigl(K(y_{i_j})K(y_{i_k})\bigr)^{\frac{n-2}{4}}},
\forall\,\,k\neq j,\,\,1\leq k,\,j\leq s. .\end{aligned}
\end{equation}
Let $\rho(\tau_s)$ be the least eigenvalue of $M(\tau_{s}),\,\forall\,\,s\in \N^\ast.$\\
$\mathbf{(A_2)}$ Assume that $\rho(\tau_{s})\neq 0$ for each
distinct points $y_{i_1},\dots,y_{i_s}\in \mathcal{K}.$\\\\We denote
by\begin{align*}\mathcal{C}_{\infty}:=\{\,&\tau_{p}:=(y_{i_1},\,\dots,\,y_{i_p}),\,p\geq
1,\,\,\mathrm{s.t}\,\,y_{i_j}\in \mathcal{K},\,\,
\forall\,\,j=1,\dots,p,\\&{}y_{i_j}\neq y_{i_k}, \forall\,\,j\neq
k,\,\,\mathrm{and\,\,}\rho(\tau_{p})> 0\,\}\end{align*}and define an
index$$i:\,\mathcal{C}_{\infty}\rightarrow
\Z,\,\,(y_{i_1},\,\dots,\,y_{i_p})\mapsto
i(y_{i_1},\,\dots,\,y_{i_p}):=p-1+\sum_{j=1}^pn-\widetilde{i}(y_{i_j}),$$where
$\widetilde{i}(y_{i_j}):=\sharp \{\,1\leq k\leq n,\,\mathrm{such
\,that\,\,}b_k(y_{i_j})< 0\,\}.$\\\\
The main result of this paper is the
following\begin{theorem}\label{t1}{\it Let $\Omega \subset
\R^n,\,\,n\geq 4,$ be a smooth bounded domain, and $\,0< K \in
\mathcal{C}^1\bigl(\bar{\Omega}\bigr)$ satisfying the assumptions
$\mathbf{(f)}_\beta,\,\mathbf{(A_1)}\,\,\mathrm{and\,\,}\mathbf{(A_2)}.$\\If$$\sum\limits_{\tau_p
\,\in \, \mathcal{C}_\infty}\,\,(-1)^{i(\tau_p)} \neq 1,$$then the
problem (\ref{problem}) has a solution.}
\end{theorem}

Our argument uses a careful analysis of the lack of compactness of
the Euler Lagrange functional $J$ associated to the problem
(\ref{problem}). Namely, we study the noncompact orbits of the
gradient flow of $J$ the so called $critical \,points\, at\,
infinity$ following the terminology of A. Bahri \cite{5}. With
respect to the closed case, new difficulties here arise. For
example, in \cite{A}, T. Aubin showed that if
$\displaystyle\frac{\partial K}{\partial \nu}\geq 0
\,\,\mathrm{on\,\,} \partial \Omega,$ then we have a possibility of
any concentration points on $\partial \Omega$  of a sequence of
subcritical solutions ( see proposition 6.44 of \cite{A} ). Using
the assumption $\mathbf{(A_1)}$ we can prove in our situation that
the boundary does not make any contribution to the existence of a
critical point at infinity. The critical points at infinity of our
problem (\ref{problem}) can be treated as usual critical points once
a Morse lemma at infinity is performed from which we can derive just
as in the classical Morse theory the difference of topology induced
by these noncompact orbits and compute their Morse index. Such a
Morse lemma at infinity is obtained through the construction of a
suitable pseudo-gradient for which the Palais-Smale condition is
satisfied along the decreasing flow lines, as long as these flow
lines do not enter the neighborhood of a finite number of critical
points $y_{i_j},\,j=1,\dots,p,$ of $K$ such that
$(y_{i_1},\,\dots,\,y_{i_p})\in \mathcal{C}_\infty.$\\

The remainder of the paper is organized as follows. In the second
section, we set up the variational structure and we recall some well
known facts. In section three, we characterize the critical points
at infinity of our problem. Section four is devoted to the proof of
the main result.

\section{Variational structure
and lack of compactness}
\def\theequation{2.\arabic{equation}}\makeatother
\setcounter{equation}{0} Our  problem (\ref{problem}) enjoys a
variational structure. Indeed, solutions of (\ref{problem})
correspond to positive critical points of the functional\be
 I(u)=\frac{1}{2}\int_{\Omega}  |\nabla u|^2 - \frac{n-2}{2n}
     \int_{\Omega} K|u|^{\frac{2n}{n-2}}
\ee defined on $H_0^1(\Omega)$.\\Let $ \Sigma:=\bigl\{\,u\in
 H^1_0(\Omega)\,,\  \mathrm{s.t.}\  \|u\|^2=\int_{\Omega}  |\nabla
u|^2=1\
  \bigr\},\,\,\Sigma^+:=\big\{\, u \in \Sigma\,,\
  u\geq 0\,\bigr\}$. Instead of working with the functional $I$ defined above, it is more convenient here to work with the
functional\be
 J(u)=\frac{\int_{\Omega}  |\nabla u|^2}
    {\Bigl( \int_{\Omega}K|u|^{\frac{2n}{n-2}}
    \Bigr)^{\frac{n-2}{n}}}
\ee defined on $\Sigma$. One can easily verify that if $u$ is a
critical point of $J$ on $\Sigma^+$, then $J(u)^\frac{n}{4}u$ is a
solution of (\ref{problem}).\\The variational viewpoint is delicate
to be studied, because the functional $J$ does not satisfy the
Palais-Smale condition ($\mathbf{(P-S)}$ for short). Which means
that there exist sequences along which $J$ is bounded, its gradient
goes to zero and which is not convergent. The analysis of the
sequences failing $\mathbf{(P-S)}$ condition can be realized
following the ideas introduced in \cite{15} and \cite{48}. For $a\in
\Omega,\,\lambda>0$, let \be
 \delta_{a,\lambda}(x)=c_n
  \Bigl(\frac{\lambda}{1+\lambda^2|x-a|^2}\Bigr)^\frac{n-2}{2},
\ee where $c_n$ is a positive constant chosen such that
 $\delta_{a,\lambda}$ is the family of solutions of the following problem
\be-\Delta u = |u|^{\frac{4}{n-2}}u ,\,\,u>0
 \,\,\mbox{in}\,
 \mathbb{R}^n.
\ee Let $P$ be the projection from $H^1(\Omega)$ on to
 $H_0^1(\Omega)$; that is, $u:= Pf$ is the unique solution of
\be \Delta u = \Delta f \,\,\qquad \mbox{in} \,\,\Omega\,,\,\qquad
\qquad u=0
 \,\,\qquad \mbox{on}
 \,\,\partial\Omega.
\ee We define now the set of potential critical points at infinity
associated to the functional $J.$ Let, for
$\e>0,\,p\in\mathbb{N}^\ast,$
\begin{align*} V(p,\e) = \bigl\{u\in\textstyle\Sigma^+ &\,  \mathrm{\, s.t\,
}
 \,  \exists \,\,  a_i\in {\Omega},\, \lambda_i> \frac{1}{\e},\,
\alpha_i>0\,  \mathrm{\, for\, } \, 1\leqslant i\leqslant p,\,
\mathrm{with}
\\
{} &\ \mathrm{ }\quad
\|u-\textstyle\sum_{i=1}^p\alpha_iP\delta_{a_i,\lambda_i}
\|<\e,\,\e_{ij}<\e,\,\,\forall\, i\neq j,
\\
{}&\quad\,\, \lambda_i{d_i}>\frac{1}{\e},\,
\Big|\textstyle\frac{\alpha_i^{\frac{4}{n-2}}K(a_i)}
{\alpha_j^{\frac{4}{n-2}}K(a_j)} -1
\Big|\!<\e\,\,\,\forall\,i,j=1,\dots,p\,\bigr\},
\end{align*}where $d_i=\mathrm{d}(a_i,\,\partial{\Omega}) \mathrm{\,\, and
\,\,}
 \e_{ij}= \bigl(
 \textstyle\frac{\lambda_i}{\lambda_j}+\frac{\lambda_j}{\lambda_i}
  +{\lambda_i}{\lambda_j}|a_i-a_j|^2
  \bigr)^{-\frac{n-2}{2}}.$\\\\If $u$ is a function in $ V(p,\e),$ one can find an optimal representation of $u$
 following the ideas introduced in \cite{6}
and \cite{10}, namely we have
\begin{proposition}\label{p2.1} For any $p\in\mathbb{N}^\ast$, there is $\e_p>0$
such that if $\e<\e_p$ and $u\in  V(p,\e)$, then the following
minimization problem \be \min\bigl\{\,
\|\,u-\textstyle\sum_{i=1}^p\alpha_iP\delta_{a_i,\lambda_i}\,\|,\,
\alpha_i>0,\,\lambda_i>0,\,a_i\in {\Omega}\,\bigr\}\ee has a unique
solution $(\bar{\alpha},\,\bar{a},\,
    \bar{\lambda})\,$(up to permutation) . Thus we can write $u$ uniquely as follows (we drop the bar):
\be u=\textstyle\sum_{i=1}^p\alpha_iP\delta_{a_i,\lambda_i}+v ,\ee
where $v$ satisfies
\begin{equation*}\label{c01}
(V_0):\,\langle v,\,\phi_i\rangle=0,\,\mathrm{for}\,\,
i=1,\dots,p,\, \mathrm{and} \,\,\phi_i=P\delta_i,\,\frac{\partial
P\delta_i}{\partial \lambda_i},\, \frac{\partial P\delta_i}{\partial
a_i}\end{equation*}
\end{proposition}
Here, $P\delta_i:=P\delta_{a_i,\lambda_i},$ and $\langle \,,\,
\rangle$ denotes the scalar product defined on $ H^1_0(\Omega)$
by$$\langle u,\,v\rangle=\int_\Omega\nabla u \nabla v.$$In the next
we will say that $v\in (V_0)$ if $v$ satisfies $(V_0).$

The failure of the $\mathbf{(P-S)}$ condition can be described
following the ideas developed in \cite{15}, \cite{48} and \cite{52}.
Such a description is by now standard and reads as follows: let
$\partial J$ be the gradient of $J$.
\begin{proposition}\label{p01}{\it
Let $(u_j)_j\subset \Sigma^+$ be a  sequence such that $\partial J$
tends to zero and $ J(u_j)$ is bounded. Then there exists an integer
$p\in\mathbb{N}^\ast$, a sequence $\e_j>0$, $\e_j\rightarrow
 0$, and an extracted subsequence of $u_j$'s, again denoted by
 $u_j$, such that $u_j\in V(p,\e_j)$}.
\end{proposition}

Now arguing as in \cite{6}, we have the following Morse lemma which
gets rid of the $v-$contribution and shows that it can be neglected
with respect to the concentration phenomenon
\begin{proposition}\label{p02}There is a $\mathcal{C}^1-$map which to each
$(\alpha_i,\,a_i,\,\lambda_i)$ such that
$\sum_{i=1}^p\alpha_iP\delta_{a_i,\lambda_i} \in V(p, \e)$
associates $\bar{v}=\bar{v}(\alpha_i,\,a_i,\,\lambda_i)$ such that
$\bar{v}$ is unique and satisfies$$J\bigl(
\sum_{i=1}^p\alpha_iP\delta_{a_i,\lambda_i}+\bar{v}
\bigr)=\mathrm{min}_{v\in (V_0)}\Bigr\{J\bigl(
\sum_{i=1}^p\alpha_iP\delta_{a_i,\lambda_i}+v \bigr)\Bigl\}.$$
Moreover, there exists a change of variables $\,v-\bar{v}\longmapsto
V\,$ such that $J$ reads in $V(p, \e)$
as$$J(\sum_{i=1}^p\alpha_iP\delta_{a_i,\lambda_i}+v)=J\bigl(
\sum_{i=1}^p\alpha_iP\delta_{a_i,\lambda_i}+\bar{v}
\bigr)+\|V\|^2.$$
\end{proposition}
The following proposition gives precise estimate of $\bar{v}.$
\begin{proposition}\label{p03}Let $\  u=\textstyle\sum_{i=1}^p\alpha_iP\delta_{a_i,\lambda_i}\in
V(p,\e),$ and let $\bar{v}$ be defined in proposition \ref{p02}. One
has the following estimate: there exists $c> 0$ independent of $u$
such that the following
holds\begin{equation}\label{v00}\begin{aligned}\|\bar{v}\|=&O\Bigl(\sum_{i=1}^p\frac{|\nabla
K(a_i)|}{\lambda_i}+\frac{1}{\lambda_i^{\beta_i}}+\frac{(\log\lambda_i)^{\frac{n+2}{2n}}}{\lambda_i^{\frac{n+2}{2}}}\Bigr)\\&{}+
\left\{
\begin{array}{ll}
O\Bigl(\sum_{i\neq
j}\e_{ij}\bigl(\log\e_{ij}^{-1}\bigr)^{\frac{n-2}{n}}+\sum_{i}
\frac{1}{(\lambda_i\,d_i)^{n-2}}\,\Bigr), & \,\,\mathrm{if}\,\,  n<6  \\
O\Bigl(\sum_{i\neq
j}\e_{ij}^{\frac{n+2}{2(n-2)}}\bigl(\log\e_{ij}^{-1}\bigr)^{\frac{n+2}{2n}}+\sum_{i}
\frac{(\log\lambda_i\,d_i)^{\frac{n+2}{2n}}}
{(\lambda_i\,d_i)^{\frac{n+2}{2}}}\,\Bigr)\quad ,
&\,\,\mathrm{if}\,\, n\geqslant 6.
\end{array}
\right.
\end{aligned}
\end{equation}
\end{proposition}
\pf Following the proof of proposition 5.3 of \cite{5}, it remains
to estimate $\int_\Omega K(x)\delta_i^{\frac{n+2}{n-2}}v\,dx,$ and
then we need to prove the following claim:
\begin{equation}\label{v0}
\int_\Omega
K(x)\delta_i^{\frac{n+2}{n-2}}v\,dx=O\Bigl(\|v\|\sum_{i=1}^p\,\frac{|\nabla
K(a_i)|}{\lambda_i}+\frac{1}{\lambda_i^{\beta_i}}+\frac{(\log\lambda_i)^{\frac{n+2}{2n}}}{\lambda_i^{\frac{n+2}{2}}}\Bigr).
\end{equation}
For this, we distinguish two cases: let $\rho> 0$ a small positive
constant such that the condition $\mathbf{(f)}_\beta$ holds in
$B(y,\,4\rho),\,\,\forall\,\,y\in \mathcal{K}.$\\\\
{\bf Case 1.} If $a_i\not\in \cup_{y\in \mathcal{K}}B(y,\,\rho).$
Let $\mu> 0$ a positive constant small enough. Using the fact that
$v\in (V_0),$ we obtain
\begin{equation}\label{v1}
\begin{aligned}\int_\Omega
K(x)\delta_i^{\frac{n+2}{n-2}}v\,dx&=\int_{B(a_i,\,\mu)}
\bigl(K(x)-K(a_i)\bigr)\delta_i^{\frac{n+2}{n-2}}v\,dx+O\Bigl(\frac{\|v\|}{\lambda_i^n}\Bigr)\\&{}\leq
c\frac{\|v\|}{\lambda_i}
\\&{}\leq c\|v\|\cdot \frac{|\nabla
K(a_i)|}{\lambda_i},\,\,\mathrm{since\,\,}\nabla K(x)\neq
0,\,\,\forall\,\,x\not\in \cup_{y\in \mathcal{K}}B(y,\,\rho).
\end{aligned}
\end{equation}
{\bf Case 2.} If $a_i \in B(y_i,\,\rho),\,\,y_i\in \mathcal{K}.$ Let
$C> 0$ a
positive constant large enough.\\
If $\lambda_i|a_i-y_i|\leq C,$ let $B_i:=B\bigl(a_i,\,\rho\bigr),$
then, by using the condition $\mathbf{(f)}_{\beta_i}$ in $B_i$ and
the fact that $v\in (V_0),$ we obtain\begin{equation}\label{v2}
\begin{aligned}
\int_\Omega
K(x)\delta_i^{\frac{n+2}{n-2}}v\,dx&=O\Bigl(\int_{B_i}\,\bigl(|x-a_i|^{\beta_i}+|y_i-a_i|^{\beta_i}\bigr)
\delta_i^{\frac{n+2}{n-2}}v\,dx\Bigr)+O\Bigl(\frac{\|v\|}{\lambda_i^{\frac{n+2}{2}}}\Bigr)
\\&{}=O\Bigl(\frac{\|v\|}{\lambda_i^{\beta_i}}+\|v\|\cdot\frac{(\log\lambda_i)^{\frac{n+2}{2n}}}{\lambda_i^{\frac{n+2}{2}}}\Bigr).
\end{aligned}\end{equation}
If $\lambda_i|a_i-y_i|\geq C,$ let
$B_i:=B\bigl(a_i,\,\frac{|a_i-y_i|}{2}\bigr)$ and
$B_{y_i}:=B\bigl(y_i,\,2\rho\bigr),$ then, by using the condition
$\mathbf{(f)}_{\beta_i}$ in $B_{y_i}$ and the fact that $v\in
(V_0),$ we obtain
\begin{align*}
\int_\Omega
K(x)\delta_i^{\frac{n+2}{n-2}}v\,dx&=\int_{B_i}\,\bigl[K(x)-K(a_i)\bigr]
\delta_i^{\frac{n+2}{n-2}}v\,dx\\&{}+\int_{B_{y_i}\setminus
B_i}\,\bigl[K(x)-K(a_i)\bigr]
\delta_i^{\frac{n+2}{n-2}}v\,dx\,+O\bigl(\frac{\|v\|}{\lambda_i^{\frac{n+2}{2}}}\bigr)
\\&{}=O\Bigl(\mathrm{sup}_{B_i}|\nabla K(x)|\cdot\int_{B_i}\,|x-a_i|\delta_i^{\frac{n+2}{n-2}}|v|\,dx\\&{}+
\int_{B_{y_i}\setminus
B_i}\,\bigl(|a_i-y_i|^{\beta_i}+|x-a_i|^{\beta_i}\bigr)\delta_i^{\frac{n+2}{n-2}}|v|\Bigr)
+O\bigl(\frac{\|v\|}{\lambda_i^{\frac{n+2}{2}}}\bigr).\end{align*}
Since $|\nabla K(x)|\sim |x-y_i|^{\beta_i-1}\,\,\mathrm{in\,\,}
B_i,$ then we obtain\begin{equation}\label{v3}\int_\Omega
K(x)\delta_i^{\frac{n+2}{n-2}}v\,dx=\|v\|\cdot O\bigl(\frac{|\nabla
K(a_i)|}{\lambda_i}+\frac{1}{\lambda_i^{\beta_i}}+\frac{(\log\lambda_i)^{\frac{n+2}{2n}}}{\lambda_i^{\frac{n+2}{2}}}\bigr).\end{equation}
Combining (\ref{v1}), (\ref{v2}) and (\ref{v3}), the claim
(\ref{v0}) follows.\\

Following A. Bahri \cite{5}, we introduce the following definition:
\begin{definition} A critical point at infinity of $J$ in $\Sigma^+$ is a limit of a flow line $u(s)$ of the equation
$$\left\lbrace
\begin{aligned}
\frac{\partial u}{\partial s}&=  -\partial J(u)\\
 u(0)&=u_0\in \Sigma^+
\end{aligned}
\right.
$$such that  $u(s)$ remains in $V(p,\e(s))$,  for $s\geq
s_0$.
\end{definition}
Here, $\e(s)$ is  some function tending to zero when $s\rightarrow
  +\infty$. Using proposition \ref{p2.1}, $u(s)$ can be written as
$$u(s)=\sum_{i=1}^p\alpha_i(s)P\delta_{a_i(s),\lambda_i(s)}+v(s).$$
Denoting by $\,a_i:=\lim a_i(s)\,\mathrm{\,and\,} \,\alpha_i:=\lim
\alpha_i(s)\,$, we denote by$$ (a_1,\dots,\,a_p)_\infty
\,\mathrm{\,or\,} \sum_{i=1}^p\alpha_iP\delta_{a_i ,\,\infty}$$ such
a critical point at infinity.

For such a critical point at infinity there are associated stable
and unstable manifolds. These manifolds can be easily described once
a Morse type reduction is performed ( see \cite{6}, pages 356-357).
\section{Characterization of the critical points at infinity }
\def\theequation{3.\arabic{equation}}\makeatother
\setcounter{equation}{0}

This section is devoted to the characterization of the critical
points at infinity, associated to the problem (\ref{problem}), in
$V(p,\e),\,\,p\geq 1.$ This characterization is obtained through the
construction of a suitable pseudo-gradient at infinity in $V(p,\e).$
The construction is based on very delicate expansion of the gradient
of the associated Euler-Lagrange functional $J$ near infinity. In
the second subsection, we will characterize the critical points at
infinity in $V(p,\e),\,\,p\geq 1.$

Using proposition  \ref{p02}, we can write, for
$u=\sum_{i=1}^p\alpha_iP\delta_{a_i\lambda_i} +v \in V(p,\e)$,$$J(u)
= J(\sum_{i=1}^p\alpha_iP\delta_{a_i\lambda_i} +\bar{v}) +
\|V\|^2.$$In the $V-$variable, we define a pseudo-gradient by
setting$$\frac{\partial V}{\partial s} = - \mu V$$ where $\mu$ is a
very large constant. Then, at $s=1$, $V(1) = e^{-\mu }V(0)$ will be
very small, as we wish. This shows that, in order to define our
deformation, we can work as if $V$ was zero. The deformation will
extend immediately, with the same properties, to a neighborhood of
zero in the $V$ variable. Therefore we need to define a
pseudo-gradient in $\{\, \sum_{i=1}^p\alpha_iP\delta_{a_i\lambda_i}
+\overline{v}\in V(p,\e)\}.$\subsection{Expansion of the gradient of
the functional }
\begin{proposition}\label{p1} \,\,Let $n\geq 4$. For $\e$ small
enough and $u = \sum_{i=1}^p\alpha_iP\delta_i \in V(p,\e)$, we have
the following expansion:
\begin{description}
\item[(i)]
\begin{align*}
\bigl\langle \partial J(u), \lambda_i \frac{\partial P \delta_i}{
\partial\lambda_i}\bigr\rangle_{H_0^1} = 2c_2J(u)\Bigl[&- \frac{n-2}{2}\alpha_i\frac{H(a_i,
a_i)}{\lambda_i^{n-2}}- \sum_{j\neq
i}\alpha_j(\lambda_i\frac{\partial \e_{ij}}{\partial\lambda_i} +
\frac{n-2}{2}\frac{H(a_i,
a_j)}{(\lambda_i\lambda_j)^{\frac{n-2}{2}}})\Bigr]\times
\\&{}\bigl(1+o(1)\bigr)+ o\Bigl(\frac{1}{\lambda_i}+\sum_{i \neq
j}\e_{ij} +\sum_{k \neq j}\e_{kj}^{\frac{n-1}{n-2}}  +
\sum_{k=1}^{p}\frac{1}{(\lambda_kd_k)^{n-1}}\Bigr),\end{align*}
where
$c_2=c_n^{\frac{2n}{n-2}}\int_{\R^n}\frac{1}{(1+|x|^2)^{\frac{n+2}{2}}}.$
\item[(ii)] If $a_i\in
B\bigl(y_{j_i},\,\rho\bigr),\,\,\mathrm{with\,\,}y_{j_i} \in
\mathcal{K},$ and $\rho$ is a positive constant small enough so that
$\mathbf{(f)}_\beta$ holds in $B\bigl(y_{j_i},\,4\rho\bigr),$ we
have
\begin{align*}
\bigl\langle \partial J(u), \lambda_i \frac{\partial P \delta_i}{
\partial\lambda_i}\bigr\rangle_{H_0^1} = 2c_2J(u)\Bigl[&
 - \frac{n-2}{2}\alpha_i\frac{H(a_i,
a_i)}{\lambda_i^{n-2}}- \sum_{j\neq
i}\alpha_j(\lambda_i\frac{\partial \e_{ij}}{\partial\lambda_i} +
\frac{n-2}{2}\frac{H(a_i,
a_j)}{(\lambda_i\lambda_j)^{\frac{n-2}{2}}})\Bigr]\times
\\&{}\bigl(1+o(1)\bigr)+o\Bigl(\sum_{i \neq j}\e_{ij} +\sum_{k \neq
j}\e_{kj}^{\frac{n-1}{n-2}}  +
\sum_{k=1}^{p}\frac{1}{(\lambda_kd_k)^{n-1}}\Bigr)\\&{}
+\bigl(\mathrm{if\,\,}\lambda_i|a_i-y_{j_i}|\geq
C\bigr)o\bigl(\frac{\bigl|\nabla
K(a_i)\bigr|}{\lambda_i}\bigr),\end{align*}where $C$ is a positive
constant large enough.\end{description}
\end{proposition}
\pf Claim ({\bf i}) is immediate from \cite{14}. Concerning claim
({\bf ii}), regarding the estimates used to prove claim ({\bf i}),
we need to estimate the quantity $\int_\Omega
K(x)\delta_i^{\frac{n+2}{n-2}}\lambda_i \frac{\partial  \delta_i}{
\partial\lambda_i}\,dx.$ Let $C$ a positive constant large
enough.\\If $\lambda_i|a_i-z_{j_i}|\leq C,$ let
$B_i:=B\bigl(a_i,\,\rho\bigr),$ then
\begin{align*}
\int_\Omega K(x)\delta_i^{\frac{n+2}{n-2}}\lambda_i \frac{\partial
\delta_i}{
\partial\lambda_i}\,dx&=\int_{B_i}\,\bigl[K(x)-K(z_{j_i})\bigr]
\delta_i^{\frac{n+2}{n-2}}\lambda_i \frac{\partial  \delta_i}{
\partial\lambda_i}\,dx\,+o\bigl(\frac{1}{\lambda_i^\beta}\bigr)
\\&{}=\frac{n-2}{2}\bigl(\frac{1}{\lambda_i^\beta}\bigr)\sum_{j=1}^n
b_j\int_{B_i}|y_j+\lambda_i(a_i-z_{j_i})_j|^\beta\frac{|y_1|^\beta\bigl(1-|y|^2\bigr)}{\bigl(|y|^2+1\bigr)^{n+1}}\,dy\,
+o\bigl(\frac{1}{\lambda_i^\beta}\bigr)\\
&{}=o\bigl(\frac{1}{\lambda_i^{n-2}}\bigr),\,\,\mathrm{since\,\,}\beta>
n-2.\end{align*}If $\lambda_i|a_i-y_{j_i}|\geq C,$ let $M> 0$ a
positive constant large enough,
$B_{i,\,k}:=B\bigl(a_i,\,\frac{|(a_i-y_{j_i})_k|}{2M}\bigr),$\\$\forall\,\,k=1,\,\dots,\,n,$
and $B_{y_{j_i}}:=B\bigl(y_{j_i},\,2\rho\bigr),$ then
\begin{equation}\label{E0}
\int_\Omega K(x)\delta_i^{\frac{n+2}{n-2}}\lambda_i \frac{\partial
\delta_i}{
\partial\lambda_i}\,dx=\int_{B_{y_{j_i}}}\,\bigl[K(x)-K(a_i)\bigr]
\delta_i^{\frac{n+2}{n-2}}\lambda_i \frac{\partial  \delta_i}{
\partial\lambda_i}\,dx+o\bigl(\frac{1}{\lambda_i^{n-2}}\bigr).\end{equation} By the
condition $\mathbf{(f)}_\beta,$ we get
\begin{equation}\label{E1}\begin{aligned}
&\quad\int_{B_{y_{j_i}}}\,\bigl[K(x)-K(a_i)\bigr]
\delta_i^{\frac{n+2}{n-2}}\lambda_i \frac{\partial  \delta_i}{
\partial\lambda_i}\,dx\\&{}=\sum_{k=1}^n\,b_k\,\int_{B_{y_{j_i}}}\bigl|(a_i-y_{j_i})_k-(a_i-x)_k\bigr|^\beta
\delta_i^{\frac{n+2}{n-2}}\lambda_i \frac{\partial  \delta_i}{
\partial\lambda_i}\,dx-\bigl|(a_i-y_{j_i})_k\bigr|^\beta\int_{B_{y_{j_i}}}
\delta_i^{\frac{n+2}{n-2}}\lambda_i \frac{\partial  \delta_i}{
\partial\lambda_i}\,dx\\&{}\quad+o\bigl(\frac{1}{\lambda_i^{n-2}}+\frac{|a_i-z_{j_i}|^{\beta-1}}{\lambda_i}\bigr)
 \end{aligned}\end{equation}
Observe that
\begin{equation}\label{E2}\begin{aligned}
&\quad\int_{B_{y_{j_i}}\setminus B_{i,\,k}
}\bigl|(a_i-y_{j_i})_k-(a_i-x)_k\bigr|^\beta
\delta_i^{\frac{n+2}{n-2}}\lambda_i \frac{\partial  \delta_i}{
\partial\lambda_i}\,dx\\&{}=O\bigl(\int_{B_{y_{j_i}}\setminus B_{i,\,k}
}(\bigl|(a_i-y_{j_i})_k\bigr|^\beta+\bigl|(a_i-x)_k\bigr|^\beta)
\delta_i^{\frac{2n}{n-2}}\,dx\bigr)=o\bigl(\frac{1}{\lambda_i^{n-2}}\bigr).
\end{aligned}\end{equation}
However, by elementary calculation, we obtain
\begin{equation}\label{E3}\begin{aligned}
&\quad\int_{B_{i,\,k}}\bigl|(a_i-y_{j_i})_k-(a_i-x)_k\bigr|^\beta
\delta_i^{\frac{n+2}{n-2}}\lambda_i \frac{\partial  \delta_i}{
\partial\lambda_i}\,dx-\bigl|(a_i-y_{j_i})_k\bigr|^\beta\int_{B_{i,\,k}}
\delta_i^{\frac{n+2}{n-2}}\lambda_i \frac{\partial  \delta_i}{
\partial\lambda_i}\,dx\\&{}=o\bigl(\frac{1}{\lambda_i^{n-2}}+\frac{|a_i-y_{j_i}|^{\beta-1}}{\lambda_i}\bigr).
\end{aligned}\end{equation}
Combining (\ref{E0}), (\ref{E1}), (\ref{E2}) and (\ref{E3}), we get
$$\int_\Omega
K(x)\delta_i^{\frac{n+2}{n-2}}\lambda_i \frac{\partial \delta_i}{
\partial\lambda_i}\,dx=o\bigl(\frac{1}{\lambda_i^{n-2}}+\frac{|a_i-y_{j_i}|^{\beta-1}}{\lambda_i}\bigr).$$
We remark from the condition $\mathbf{(f)}_\beta$ that, for $\rho$
small enough, $\bigl|\nabla K(a_i)\bigr|\sim
\bigl|(a_i-y_{j_i})\bigr|^{\beta-1}.$ Then we can write$$\int_\Omega
K(x)\delta_i^{\frac{n+2}{n-2}}\lambda_i \frac{\partial \delta_i}{
\partial\lambda_i}\,dx=o\bigl(\frac{1}{\lambda_i^{n-2}}+\frac{\bigl|\nabla K(a_i)\bigr|}{\lambda_i}\bigr).$$
This finishes the proof of claim ({\bf ii}).
\begin{proposition}\label{p2} Let $n\geq 4$. For $\e$ small
enough and $u = \sum_{i=1}^p\alpha_iP\delta_i \in V(p,\e)$, we have
\begin{description}
\item[(i)]
\begin{align*}&\bigl\langle \partial J(u),\frac{1}{\lambda_i}
\frac{\partial P \delta_i}{\partial a_i}\bigr\rangle_{H_0^1} =
2J(u)\Bigl[
 -\alpha_i^{\frac{n+2}{n-2}}\frac{n-2}{n}c_4J^{\frac{n}{n-2}} \frac{\nabla
K(a_i)}{\lambda_i}+
\frac{\alpha_i}{\lambda_i^{n-1}}c_2\frac{\partial H(a_i,
a_i)}{\partial a_i}\\&\phantom{( \partial J(u),\frac{1}{\lambda_i}
\frac{\partial P \delta_i}{\partial a_i})_{H_0^1}= J(u)\Bigl[ {}}-
c_2 \sum_{j\neq i}\alpha_j(\frac{1}{\lambda_i}\frac{\partial
\e_{ij}}{\partial a_i} -
\frac{1}{(\lambda_i\lambda_j)^{\frac{n-2}{2}}\lambda_i}
\frac{\partial H(a_i, a_j)}{\partial
a_i})\Bigr](1+o(1))\\&\phantom{(
\partial J(u),\frac{1}{\lambda_i} \frac{\partial P
\delta_i}{\partial a_i})_{H_0^1}= J(u)\Bigl[ {}}
+o\bigl(\frac{1}{\lambda_i}\sum_{k \neq
j}\e_{kj}^{\frac{n-1}{n-2}}+\sum_{k=1}^{p}\frac{1}{(\lambda_kd_k)^{n-1}}\bigr)+
O\Bigl(\lambda_j|a_i-a_j|\e_{ij}^{\frac{n+1}{n-2}}\Bigr),\end{align*}where
$c_4=c_n^{\frac{2n}{n-2}}\int_{\R^n}\frac{|x|^2}{(1+|x|^2)^{n+1}}.$
\item[(ii)] If $\,a_i\in
B\bigl(y_{j_i},\,\rho\bigr),\,\,\mathrm{with\,\,}y_{j_i} \in
\mathcal{K},$ and $\rho$ is a positive constant small enough so that
$\mathbf{(f)}_\beta$ holds in $B\bigl(y_{j_i},\,4\rho\bigr),$ then
the above estimate can be improved. Let $C$ a positive constant
large enough. If $\lambda_i|(a_i-z_{j_i})_k|\geq C,$ we get
\begin{align*}\bigl\langle \partial J(u),\frac{1}{\lambda_i}
\frac{\partial P \delta_i}{\partial (a_i)_k}\bigr\rangle_{H_0^1} &=
-2(n-2)J^{\frac{2n-2}{n-2}}\alpha^{\frac{n+2}{n-2}}\mathrm{sg}\bigl[(a_i-z_{j_i})_k\bigr]
\frac{\bigl|(a_i-z_{j_i})_k\bigr|^{\beta-1}}{\lambda_i}b_kc_5
\\&{}+o\bigl(\frac{1}{\lambda_i^{n-2}}+\frac{\bigl|\nabla
K(a_i)\bigr|}{\lambda_i}+\sum_{j \neq i}\e_{ij}+\sum_{k \neq
j}\e_{kj}^{\frac{n-1}{n-2}}+\sum_{k=1}^{p}\frac{1}{(\lambda_kd_k)^{n-1}}\bigr)
\Bigr)\\&{}+ O\Bigl( \sum_{j\neq
i}\bigl|\frac{1}{\lambda_i}\frac{\partial \e_{ij}}{\partial
a_i}\bigr|.\end{align*} Here,
$c_5=c_n^{\frac{2n}{n-2}}\frac{\beta}{n}\int_{\R^n}\,\frac{|x|^2}{\bigl(1+|x|^2
\bigr)^{n+1}}\,dx,$ $k\in \{\,1,\,\dots,\,n\,\},$ and $(a_i)_k$
denotes the $k^{t\,h}$ component of $a_i$ in some local coordinates
system.\end{description}
\end{proposition}
\pf Claim ({\bf i}) is immediate from \cite{14}. Concerning claim
({\bf ii}), arguing as in the proof of proposition \ref{p1}, claim
({\bf ii}) is proved under the following estimate: let $C$ a
positive constant large enough and $a_i\in
B\bigl(y_{j_i},\,\rho\bigr),$ where $\rho$ is a positive constant
small enough so that $\mathbf{(f)}_\beta$ holds in
$B\bigl(y_{j_i},\,4\rho\bigr).$ If $\lambda_i|(a_i-z_{j_i})_k|\geq
C,$ then we have\begin{align*} \int_\Omega
K(x)\delta_i^{\frac{n+2}{n-2}}\frac{1}{\lambda_i} \frac{\partial
\delta_i}{
\partial
(a_i)_k}\,dx&=(n-2)\mathrm{sg}\bigl[(a_i-z_{j_i})_k\bigr]\frac{\bigl|(a_i-z_{j_i})_k\bigr|^{\beta-1}}{\lambda_i}b_kc_5+
o\bigl(\frac{1}{\lambda_i^{n-2}}+\frac{\bigl|\nabla
K(a_i)\bigr|}{\lambda_i}\bigr) .\end{align*} This finishes the proof
of claim ({\bf ii}).
\subsection{Critical points at infinity }

This subsection is devoted to the characterization of the critical
points at infinity , associated to the problem (\ref{problem}), in
$V(p,\e),\,\,p\geq 1.$ This characterization is obtained through the
construction of a suitable pseudo-gradient at infinity for which the
Palais-Smale condition is satisfied along the decreasing flow lines
as long as these flow lines do not enter in the neighborhood of
finite number of critical points $y_{i_j},\,j=1,\dots,p,$ of $K$
such that $(y_{i_1},\,\dots,\,y_{i_p})\in \mathcal{C}_\infty.$ Now,
we introduce the following main result:
\begin{theorem}\label{t2}\,\,Let $n\geq4$. There exists a
pseudo-gradient $W$ so that the following holds. There is a constant
$c>0$ independent of $u = \sum_{i=1}^p\alpha_iP\delta_{a_i\lambda_i}
\in V(p,\e)$ so that:\\{\bf (i)} $$\Bigl\langle\partial J(u), W(u)
\Bigr\rangle \leq -c\Bigl(\sum_{i}^{p} [ \frac{|\nabla
K(a_i)|}{\lambda_i} + \frac{1}{\lambda_i^{n-2}} +
\frac{1}{(\lambda_id_i)^{n-1}} ] + \sum_{i \neq
j}\e_{ij}^{\frac{n-1}{n-2}} \Bigr).$${\bf (ii)}
$$\Bigl\langle\partial J(u+\overline{v}),
W(u)+\frac{\partial\overline{v}}{\partial(\alpha, a,
\lambda)}(W)\Bigr\rangle
 \leq -c\Bigl(\sum_{i}^{p} [ \frac{|\nabla K(a_i)|}{\lambda_i}
+ \frac{1}{\lambda_i^{n-2}} + \frac{1}{(\lambda_id_i)^{n-1}} ] +
\sum_{i \neq j}\e_{ij}^{\frac{n-1}{n-2}} \Bigr).$${\bf (iii)} The
minimal distance to the boundary, $d_i(t): = \mathrm{d}(a_i(t),
\partial \Omega),$ only increases if it is small
enough.
\\{\bf (iv)} $|W|$ is bounded. Furthermore, the only case where the
maximum of the $\lambda_i$'s is not bounded is when each point $a_j$
is close to a critical point $y_{i_j}$ with $y_{i_j} \neq y_{i_k},$
for each $j \neq k,$ and $(y_{i_1},\,\dots,\,y_{i_p})\in
\mathcal{C}_{\infty}.$\end{theorem}Before giving the proof of
theorem \ref{t2}, we need to state two results which deal with two
specific cases of theorem \ref{t2}. The proof of these results will
be given later. Let $d_0>0$ be a constant small enough such
that$$\frac{\partial K}{\partial \nu} (x) < -c_0,\,\,
 \forall x \in \Omega_{d_0} := \{ x\in \Omega, d(x, \partial
 \Omega) \leq 2d_0 \},$$where $c_0$ is a fixed positive constant. Then we have the following
propositions:\begin{proposition}\label{p3}In
$V_{d_0}(p,\,\e):=\{u=\sum_{i=1}^p\alpha_iP\delta_i \in
V(p,\e),\,d(a_i,
\partial \Omega)\geq d_0,\,\forall \,1\leq i\leq p\},$ there
exists a pseudo-gradient $W_1$ so that the following holds: There is
a constant $c > 0$ independent of $u \in V_{d_0}(p,\,\e)$ so
that$$\Bigl\langle\partial J(u), W_1(u)\Bigr\rangle \leq
-c\Bigl(\sum_{i}^{p}\, [ \frac{|\nabla K(a_i)|}{\lambda_i} +
\frac{1}{\lambda_i^{n-2}}
 ] + \sum_{i \neq
j}\e_{ij} \Bigr).$$
\end{proposition}\begin{proposition}\label{p4} In $V_b(p,\e):= \{ u = \sum_{i=1}^p\alpha_iP\delta_i \in
V(p,\,\e),\, d(a_i, \partial \Omega) \leq 2d_0 ,\, \forall \,\,
1\leq i\leq p \},$ there exists a pseudo-gradient $W_2$ so that the
following holds: There is a constant $c > 0$ independent of $u \in
V_b(p,\,\e)$ so that
$$\Bigl\langle\partial J(u), W_2(u)\Bigr\rangle \leq -c\Bigl(\sum_{i}^{p}\, [\frac{1}{\lambda_i} +
 \frac{1}{(\lambda_id_i)^{n-1}}] + \sum_{i \neq j}\e_{ij}^{\frac{n-1}{n-2}}\Bigr).$$
\end{proposition}
\pfn {\bf theorem \ref{t2}} We divide the set $\{ 1, ..., p\}$ into
 two sets. The first contains the indices of the points near
 the boundary $\partial \Omega$, and the second contains the indices of
 the points far away from $\partial \Omega$. Let us define

$B:= \{1\leq i\leq p \,\,\mathrm{ s.t}\, \, d_i \geq 2d_0\}.$\\

$B_1:= B \cup \{ i \not\in B\,\,\mathrm{s.t}\,\, \exists (i_1, ...,
i_r) \,\,\mathrm{with}\, \, i_1=i, \,i_r\in B\,\,\mathrm{ and}\,\,
|a_{i_{k-1}}-a_{i_{k}}| < \frac{d_0}{p},\,\, \forall \,\, k \leq r\}.$\\

$B_2:= \{ 1,\,\dots,\, p\} \backslash B_1.$\\\\
Observe that

{\bf $(O_1)$ } $ d_i:=d(a_i, \partial \Omega) \leq 2d_0 ,\,\,\forall \,\,i\in B_2.$\\

{\bf $(O_2)$ } The advantage of $B_1$ is that if $i \in B_1$ and
$j\not\in B_1,$ then $|a_i-a_j| \geq \frac{d_0}{p}.$\\Now we write
$u$ as$$u:= u_1 + u_2, \,\, u_k:= \sum_{i \in
 B_k}\alpha_iP\delta_i \,\quad (1 \leq k \leq 2).$$Observe that $u_1\in V_{d_0}(\mathrm{card}(B_1), \e)$. Then we use
the previous construction as in proposition \ref{p3} to $u_1$, which
means we apply the previous construction to the sub-pack of
functions $u := \sum_{i =1}^{\mathrm{card}(B_1)}\,\alpha_iP\delta_i$
forgetting the indices $i \not\in B_1$. Let $W_1(u_1)$ be the vector
field thus defined. The same argument can be repeated for $u_2$,
which is in $V_b(\mathrm{card}(B_2), \e)$, and we will denote by
$W_2(u_2)$ the vector field thus defined. Define $W$ as $W(u)=
W_1(u_1)+W_2(u_2).$ Thus we
have\begin{equation*}\label{}\begin{aligned}\Bigl\langle\partial
J(u), W(u)\Bigr\rangle &=\Bigl\langle\partial J(u_1), W_1(u_1)
\Bigr\rangle+\Bigl\langle\partial
J(u_2),\,W_2(u_2)\Bigr\rangle\\&{}+o\Bigl(\sum_{k=1}^{p}[\frac{1}{\lambda_k^{n-2}}
+ \frac{1}{(\lambda_kd_k)^{n-1}} ] \Bigr)+O\Bigl( \sum_{i\in
B_1,\,j\in B_2}\e_{ij}\Bigr).\end{aligned}\end{equation*}Observe
that, for $i\in B_1\,\,\mathrm{and\,\,}j\in B_2,$$$\e_{ij}=o\bigl(
\frac{1}{\lambda_i^{n-2}}+\frac{1}{\lambda_j}\bigr).$$So claim {\bf
(i)} of theorem \ref{t2} follows. Now, arguing as in appendix 2 of
\cite{6}, claim {\bf(ii)} follows from {\bf(i)} and proposition
\ref{p03}. The conditions {\bf(iii)} and {\bf(iv)} are satisfied by
the definition of the vector field $W.$ \vskip3mm\n
\pfn {\bf proposition \ref{p3}.} Let $\eta >0$ a fixed constant small enough with $|y_i-y_j|>2\eta \,\,\forall \,\,i\neq j.$\\
We divide the set $V(p,
 \e,d_0)$ into three sets:\begin{align*}\label{}
\begin{split}V_1(p,\e,d_0):=\{& u = \sum_{i=1}^p\alpha_iP\delta_i \in
V_{d_0}(p,\e)\,\,\mbox{s.t}\,\, a_i\in B(y_{j_i}, \eta),\,y_{j_i}\in \mathcal{K},\,\,\forall \,i=1,\dots,p,\\
&\text{with} \,y_{j_i}\neq y_{j_k} \,\forall\,i\neq k,\,\,
\mbox{and}\,\, \rho(y_{j_1},\dots,\,y_{j_p})> 0\}.
\end{split}\\
\begin{split}
V_2(p,\e,d_0) := \{& u = \sum_{i=1}^p\alpha_iP\delta_i \in
V_{d_0}(p,\e)\,\,\mbox{s.t}\,\, a_i\in B(y_{j_i}, \eta),\,y_{j_i}\in \mathcal{K},\,\,\forall \,i=1,\dots,p,\\
&\text{with} \,y_{j_i}\neq y_{j_k} \,\forall\,i\neq k
,\,\,\mbox{and}\,\, \rho(y_{j_1},\dots,\,y_{j_p})< 0\}.
\end{split}\\\begin{split}V_3(p,\e,d_0):= \{&u = \sum_{i=1}^p\alpha_iP\delta_i \in
V_{d_0}(p,\e)\,\,\mbox{s.t}\,\, a_i \in B(y_{j_i},
\eta),\,y_{j_i}\in \mathcal{K},\,\,\forall \,i=1,\dots,p,\,
\\&\text{and}\,\,\exists \,\,i\neq k \,\mbox{s.t} \,y_{j_i}=y_{j_k} \}.\end{split}\\
V_4(p,\e,d_0) := \{& u = \sum_{i=1}^p\alpha_iP\delta_i \in
V_{d_0}(p,\e)\,\,\text{s.t}\,\,\exists \,\, a_i\not\in \cup_{y\in
\mathcal{K}}\,B(y, \eta)\}.
\end{align*}We will define the pseudo-gradient depending on the sets $V_i(p,
\e),\,i=1-4,$ to which $u$ belongs.\begin{lemma}\label{l1}\,In
$V_2(p,\e,d_0)$, there exists a pseudo-gradient $\widetilde{W}_2$ so
that the following holds: There is a constant $c > 0$ independent of
$u \in V_2(p,\e,d_0)$ so that$$\bigl\langle
\partial J(u), \widetilde{W}_2\bigr\rangle\leq
-c\Bigl(\sum_{i=1}^p\frac{1}{\lambda_i^{n-2}}+\frac{\bigl|\nabla
K(a_i)\bigr|}{\lambda_i}+\sum_{i \neq j}\e_{ij}\Bigr).$$\end{lemma}
\pf Let $\rho$ be the least eigenvalue of $M$. Then there exists an
eigenvector $e=(e_1,\dots,\,e_p)$ associated to $\rho$ such that
$\|e\|=1$ with $e_i> 0,\,\forall\,i=1,\dots,\,p.$ Indeed, let
$e=(e_1,\dots,\,e_p)$ an eigenvector associated to $\rho,$ with
$\|e\|=1.$ By elementary calculation, we obtain, since $m_{ij}<
0,\,\forall\,i\neq j,$\begin{equation}\label{equations1}e_i>
0,\,\forall\,\,i=1,\dots,\,p,\,\,\mathrm{or\,\,}e_i<
0,\,\forall\,\,i=1,\dots,\,p.\end{equation} Let $\gamma> 0$ such
that for any $x\in B(e,\,\gamma):=\{y\in S^{p-1}\,/\,\|y-e\|\leq
\gamma\}$ we have $^Tx\cdot M\cdot x< (1/2)\rho.$ Two cases may
occur:\\$\bf Case\, 1:$ $|\Lambda|^{-1}\Lambda\in B(e,\,\gamma).$
Since, for any $i\neq j$, $|a_i-a_j|\geq c$, then
\begin{equation}\label{c1}\lambda_i\frac{\partial
\e_{ij}}{\partial
\lambda_i}=-\e_{ij}(1+o(1))=-\frac{1}{\bigl(\lambda_i\lambda_j|a_i-a_j|^2\bigr)^{\frac{n-2}{2}}}(1+o(1))\end{equation}and
\begin{equation}\label{c2}\frac{1}{\lambda_i}\frac{\partial
\e_{ij}}{\partial a_i}=o\bigl(\e_{ij}\bigr).\end{equation} We define
$\widetilde{W}_2^{1}$ by $$\widetilde{W}_2^{1} =
-\sum_{i=1}^p\alpha_i\lambda_i\frac{\partial P
\delta_i}{\partial\lambda_i}.$$From proposition \ref{p1} and
(\ref{c1}), we obtain
\begin{equation*}\label{}
\begin{aligned}
\bigl\langle \partial J(u), \widetilde{W}_2^{1}\bigr\rangle &=
-64\pi^2J(u)\Biggl[\sum_{i=1}^p\alpha_i^2\frac{H(a_i,a_i)}{\lambda_i^{n-2}}-\sum_{j\neq
i}\alpha_j\alpha_i\lambda_i\frac{G(a_i,a_j)
}{\bigl(\lambda_i\lambda_j\bigr)^{\frac{n-2}{2}}}\Biggr]\times
\bigl(1+o(1)\bigr)
\\&+\sum_{i=1}^p\,
\bigl(\mathrm{if\,\,}\lambda_i|a_i-y_{j_i}|\geq
C\,\bigr)o\Bigl(\frac{\bigl|\nabla K(a_i)\bigr|}{\lambda_i}\Bigr).
\end{aligned}
\end{equation*} Observe that, since $u\in V(p,\e)$, we have $J(u)^{\frac{n}{n-2}}\alpha_i^{\frac{4}{n-2}}K(a_i)= (1+o(1)).$ Thus we derive that
\begin{equation}\label{c3}
\begin{aligned}
\bigl\langle\partial J(u), \widetilde{W}_2^{1}\bigr\rangle & =
-c\Bigl[^T\Lambda \cdot M\cdot \Lambda
\Bigr]\bigl(1+o(1)\bigr)+\sum_{i=1}^p\bigl(\mathrm{if\,}\lambda_i|a_i-y_{j_i}|\geq
C\bigr)o\Bigl( \sum_{i=1}^p\frac{\bigl|\nabla
K(a_i)\bigr|}{\lambda_i}\Bigr)\\&{}\leq-c\bigl(\sum_{i=1}^p\frac{1}{\lambda_i^{n-2}}\bigr)
+\sum_{i=1}^p\bigl(\mathrm{if\,\,}\lambda_i|a_i-y_{j_i}|\geq
C\,\bigr)o\Bigl( \frac{\bigl|\nabla K(a_i)\bigr|}{\lambda_i}\Bigr)
\\&{}\leq-c\bigl(\sum_{i=1}^p\frac{1}{\lambda_i^{n-2}}+\sum_{i \neq
j}\e_{ij}\bigr)+\sum_{i=1}^p\bigl(\mathrm{if\,\,}\lambda_i|a_i-y_{j_i}|\geq
C\,\bigr)o\Bigl( \sum_{i=1}^p\frac{\bigl|\nabla
K(a_i)\bigr|}{\lambda_i}\Bigr),
\end{aligned}
\end{equation}
where $M$ is the matrix defined by (\ref{matrice}), and
$\Lambda:=^T(\frac{1}{\lambda_1^{\frac{n-2}{2}}},\dots,\,\frac{1}{\lambda_p^{\frac{n-2}{2}}})$.\\
$\bf Case\, 2:$ $|\Lambda|^{-1}\Lambda\not\in B(e,\,\gamma).$ In
this case, we define
$$\widetilde{W}_2^{2} =
-\sum_{i=1}^p|\Lambda|\alpha_i\beta_i\lambda_i^{\frac{n}{2}}\frac{\partial
P \delta_i}{\partial\lambda_i},$$where
$\beta_i=\frac{|\Lambda|e_i-\Lambda_i}{\|y(0)\|}-\frac{y_i(0)}{\|y(0)\|^3}(y(0),\,|\Lambda|e-\Lambda)
\,\mbox{and}\,y(t)=(1-t)\Lambda+t|\Lambda|e.$ Define
$$\Lambda(t):=y(t)/\|y(t)\|.$$ Using proposition \ref{p1}, we
derive\begin{equation}\label{equations2}
\begin{aligned}
\bigl\langle \partial J(u), \widetilde{W}_2^{2}\bigr\rangle &=
cJ(u)|\Lambda|\Biggl[\sum_{i=1}^p\alpha_i^2\beta_i\frac{H(a_i,a_i)}{\lambda_i^{\frac{n-2}{2}}}-\sum_{j\neq
i}\alpha_j\alpha_i\beta_i\frac{G(a_i,a_j)
}{\lambda_j^{\frac{n-2}{2}}}\Biggr] \\&{}+ o\Bigl(\sum_{i \neq
j}\e_{ij} +
\sum_{i=1}^p\frac{1}{\lambda_i^{n-2}}\Bigr)+\sum_{i=1}^p\bigl(\mathrm{if\,\,}\lambda_i|a_i-y_{j_i}|\geq
C\,\bigr)o\Bigl( \frac{\bigl|\nabla
K(a_i)\bigr|}{\lambda_i}\Bigr)\\&{}=\frac{c}{J(u)^{\frac{n-2}{2}}}|\Lambda|^2\Bigl[^T\Lambda^{'}(0)
\cdot M\cdot \Lambda(0)\Bigr]\\&{}+o\bigl(
\sum_{i=1}^p\frac{1}{\lambda_i^{n-2}}\bigr)+\sum_{i=1}^p\bigl(\mathrm{if\,}\lambda_i|a_i-y_{j_i}|\geq
C\,\bigr)o\Bigl( \frac{\bigl|\nabla
K(a_i)\bigr|}{\lambda_i}\Bigr),\,\mathrm{since\,}|a_i-a_j|\geq
c,\,\forall\,\,i\neq j.
\end{aligned}
\end{equation}We claim
that\begin{equation}\label{equations3}\frac{\partial}{\partial
t}(^T\Lambda(t) \cdot M\cdot \Lambda(t))<
-c,\,\mbox{for}\,\,t\,\,\mbox{near}\,\,0,\end{equation}where $c>0$
is a constant independent of $|\Lambda|^{-1}\Lambda \in
\bigl(B(e,\,\gamma)\bigr)^c.$ Indeed, we have
\begin{equation}\label{equations4}
^T\Lambda(t) \cdot M\cdot
\Lambda(t)=\rho+\frac{(1-t)^2}{\|y(t)\|^2}\bigl[^T\Lambda \cdot
M\cdot \Lambda-\rho|\Lambda|^2\bigr].
\end{equation}
Equation (\ref{equations4}) implies
\begin{equation}\label{equations5}\begin{aligned}
\frac{\partial}{\partial t}(^T\Lambda(t) \cdot M\cdot
\Lambda(t))&=\frac{2(1-t)}{\|y(t)\|^4}\bigl[^T\Lambda \cdot M\cdot
\Lambda-\rho|\Lambda|^2\bigr]\bigl[-(1-t)|\Lambda|(e,\,\Lambda)-t|\Lambda|^2\bigr]\\&{}
=\frac{2(1-t)|\Lambda|^4}{\|y(t)\|^4}\Bigl[\frac{1}{|\Lambda|^4}\bigl(^T\Lambda
\cdot M\cdot
\Lambda-\rho|\Lambda|^2\bigr)\bigl(-|\Lambda|(e,\,\Lambda)\bigr)+o(1)\Bigr]
\end{aligned}
\end{equation}
By using the observation (\ref{equations1}), we derive that there
exists $c> 0$ ( $c$ independent of $|\Lambda|^{-1}\Lambda \in
\bigl(B(e,\,\gamma)\bigr)^c$ ) such that
\begin{equation}\label{equations6}
^T\Lambda \cdot M\cdot \Lambda-\rho|\Lambda|^2\geq c |\Lambda|^2.
\end{equation}
Also, observe that
\begin{equation}\label{equations7}
|\Lambda|(e,\,\Lambda) \geq \alpha
|\Lambda|^2,\,\,\mathrm{where\,\,}\alpha:=\mathrm{inf}\{e_i,\,\,1\leq
i\leq p\}.
\end{equation}Combining (\ref{equations5}), (\ref{equations6}) and
(\ref{equations7}), the claim (\ref{equations3}) follows. Now, by
combining (\ref{equations2}) and (\ref{equations3}), we obtain
\begin{equation}\label{c33}
\begin{aligned}
\bigl\langle \partial J(u), \widetilde{W}_2^{2}\bigr\rangle &\leq
-c\bigl(\sum_{i=1}^p\frac{1}{\lambda_i^{n-2}}+\sum_{i \neq
j}\e_{ij}\bigr)+\sum_{i=1}^p\bigl(\mathrm{if\,\,}\lambda_i|a_i-y_{j_i}|\geq
C\,\bigr)o\Bigl( \frac{\bigl|\nabla K(a_i)\bigr|}{\lambda_i}\Bigr).
\end{aligned}
\end{equation}We define $\widetilde{W}_2^{3}$ as a convex combination of $\widetilde{W}_2^{1}$ and
$\widetilde{W}_2^{2}.$ Combining (\ref{c3}) and (\ref{c33}), we
obtain
\begin{equation}\label{c3"}
\begin{aligned}
\bigl\langle \partial J(u), \widetilde{W}_2^{3}\bigr\rangle \leq
-c\bigl(\sum_{i=1}^p\frac{1}{\lambda_i^{n-2}}+\sum_{i \neq
j}\e_{ij}\bigr)+\sum_{i=1}^p\bigl(\mathrm{if\,\,}\lambda_i|a_i-y_{j_i}|\geq
C\,\bigr)o\Bigl( \frac{\bigl|\nabla K(a_i)\bigr|}{\lambda_i}\Bigr).
\end{aligned}
\end{equation}
Let $\Psi$ a positive cut-off function defined by
$\Psi(t)=1,\,\,\mathrm{if\,\,}t\leq C$ and
$\Psi(t)=0,\,\,\mathrm{if\,\,}t\geq 2C,$ where $C$ is a positive
constant large enough. To make appear
$\sum_{i=1}^p\frac{\bigl|\nabla K(a_i)\bigr|}{\lambda_i},$ we
define, for each $i=1,\dots,\,p
,$$$\overline{X}_i=\sum_{k=1}^n\,\Bigl[1-\Psi\bigl(\lambda_i|(a_i-y_{j_i})_k|\bigr)\Bigr]b_k\cdot
\mathrm{sg}\bigl[(a_i-y_{j_i})_k\bigr] \frac{1}{\lambda_i}
\frac{\partial P \delta_i}{\partial (a_i)_k}.$$ Using proposition
\ref{p2} and (\ref{c2}), we obtain
\begin{equation}\label{c4}
\begin{aligned}\bigl\langle \partial J(u), \overline{X}_i\bigr\rangle
&=-c\sum_{k=1}^n\,\Bigl[1-\Psi\bigl(\lambda_i|(a_i-y_{j_i})_k|\bigr)\Bigr]b_k^2\frac{\bigl|(a_i-y_{j_i})_k\bigr|^{\beta_{j_i}-1}}{\lambda_i}
\\&\quad+o\bigl(\sum_{i=1}^p\frac{\bigl|\nabla
K(a_i)\bigr|}{\lambda_i}+\frac{1}{\lambda_i^{n-2}}\bigr)
\\&{}\leq -c\Bigl[1-\Psi\bigl(\lambda_i|(a_i-y_{j_i})_{k_i}|\bigr)\Bigr]\frac{\bigl|\nabla
K(a_i)\bigr|}{\lambda_i}+o\bigl(\sum_{i=1}^p\frac{\bigl|\nabla
K(a_i)\bigr|}{\lambda_i}+\frac{1}{\lambda_i^{n-2}}\bigr),
\end{aligned}
\end{equation}where $k_i$ denotes the index such that $|(a_i-y_{j_i})_{k_i}|=\mathrm{max}_{1\leq k\leq n}|(a_i-y_{j_i})_k|.$
Combining (\ref{c3"}) and (\ref{c4}), we obtain
\begin{equation}\label{c5}
\begin{aligned}\bigl\langle \partial J(u),\widetilde{W}_2^{3}+\sum_{i=1}^p\overline{X}_i\,\bigr\rangle\leq&
-c\Bigl(\sum_{i=1}^p\frac{1}{\lambda_i^{2}}+\sum_{i \neq
j}\e_{ij}+\bigl[1-\Psi\bigl(\lambda_i|(a_i-y_{j_i})_{k_i}|\bigr)\bigr]\frac{\bigl|\nabla
K(a_i)\bigr|}{\lambda_i}\Bigr)\\&{}+o\bigl(\sum_{i=1}^p\frac{\bigl|\nabla
K(a_i)\bigr|}{\lambda_i}\bigr).
\end{aligned}
\end{equation}
If $\Psi \leq \frac{1}{2},$ then $\frac{\bigl|\nabla
K(a_i)\bigr|}{\lambda_i}$ appears in the upper bound of (\ref{c5}).
However, if $\Psi > \frac{1}{2},$ then we have $\frac{\bigl|\nabla
K(a_i)\bigr|}{\lambda_i}=o\bigl( \frac{1}{\lambda_i^{n-2}}\bigr),$
and so we can make appear $\frac{\bigl|\nabla
K(a_i)\bigr|}{\lambda_i}$ from $\frac{1}{\lambda_i^{n-2}}.$ From
this discussion, the estimate (\ref{c5}) becomes
\begin{equation}\label{*10}
\begin{aligned}\bigl\langle \partial J(u),\widetilde{W}_2^{3}+\sum_{i=1}^p\overline{X}_i\,\bigr\rangle\leq&
-c\Bigl(\sum_{i=1}^p\frac{1}{\lambda_i^{2}}+\frac{\bigl|\nabla
K(a_i)\bigr|}{\lambda_i}+\sum_{i \neq j}\e_{ij}\Bigr).
\end{aligned}
\end{equation}The pseudo-gradient
$\widetilde{W}_2:=\widetilde{W}_2^{3}+\sum_{i=1}^p\,\overline{X}_i$
satisfies the claim of lemma \ref{l1}
\begin{lemma}\label{l2}\,In $V_1(p,\e,d_0)$, there exists a pseudo-gradient $\widetilde{W}_1$ so that
the following holds: There is a constant $c > 0$ independent of $u
\in V_1(p,\e,d_0)$ so that
$$\bigl\langle\partial J(u), \widetilde{W}_1(u)\bigr\rangle \leq -c\Bigl(\sum_{i=1}^p\frac{1}{\lambda_i^{n-2}}+\frac{\bigl|\nabla K(a_i)\bigr|}{\lambda_i}+\sum_{i \neq
j}\e_{ij}\Bigr).$$\end{lemma}\pf Let $\delta> 0$ a fixed constant
small enough, and denote, for each $\beta_{j_i},$
$\alpha_{j_i}:=\frac{n-3}{\beta_{j_i}-1}.$ We distinguish two cases:\\
${\bf Case\, 1}:$ $\mathrm{max}_{1\leq i\leq
p}\{\lambda_i^{\alpha_{j_i}}|a_i-y_{j_i}|\}\leq \delta.$ In this
case, we define$$Y_1 := \sum_{i=1}^p\alpha_i\lambda_i\frac{\partial
P \delta_i}{\partial\lambda_i}.$$ Arguing as in the proof of the
estimate (\ref{c3}), and using the fact
$\rho(y_{j_1},\dots,\,y_{j_p})> 0,$ we obtain
\begin{equation}\label{cR1}
\begin{aligned}
\bigl\langle\partial J(u), Y_1\bigr\rangle
\leq-c\bigl(\sum_{i=1}^p\frac{1}{\lambda_i^{n-2}}+\sum_{i \neq
j}\e_{ij}\bigr)+\sum_{i=1}^p\bigl(\mathrm{if\,\,}\lambda_i|a_i-y_{j_i}|\geq
C\,\bigr)o\Bigl( \sum_{i=1}^p\frac{\bigl|\nabla
K(a_i)\bigr|}{\lambda_i}\Bigr).
\end{aligned}
\end{equation}
Observe that $\frac{\bigl|\nabla
K(a_i)\bigr|}{\lambda_i}=o\bigl(\frac{1}{\lambda_i^{n-2}}\bigr),\,\forall\,1\leq
i\leq p.$ Thus, from (\ref{cR1}), we get
\begin{equation}\label{cR2}
\begin{aligned}
\bigl\langle\partial J(u), Y_1\bigr\rangle
\leq-c\bigl(\sum_{i=1}^p\frac{1}{\lambda_i^{n-2}}+\frac{\bigl|\nabla
K(a_i)\bigr|}{\lambda_i}+\sum_{i \neq j}\e_{ij}\bigr).
\end{aligned}
\end{equation}
$\bf Case\, 2:$ $\mathrm{max}_{1\leq i\leq
p}\{\lambda_i^{\alpha_{j_i}}|a_i-y_{j_i}|\}> \delta.$
Let$$i_1:=\mathrm{min}\{\,1\leq i\leq
p,\,\,\mathrm{s.t\,\,}\lambda_i^{\alpha_{j_i}}|a_i-y_{j_i}|\}>
\delta\,\}.$$Without loss of generality, we suppose $\lambda_1\leq
\dots\leq \lambda_p.$ Let $M> 0$ a fixed constant large enough. We
set\begin{equation}\label{cR3}
\begin{aligned}I&:=\{i_1\}\cup \{\,i< i_1,\,\,\mathrm{s.t\,\,}\lambda_{j-1}\geq
\frac{1}{M}\lambda_j,\,\forall\,\,i< j\leq
i_1\,\}\\&=:\{i_0,\dots,\,i_1\}.\end{aligned}
\end{equation}
We distinguish two subcases:\\ ${\bf Subcase\, 2.1}:$ $i_0> 1.$ Let
$$\widetilde{u}:=\sum_{i< i_0}\alpha_iP\delta_i.$$ $\widetilde{u}$
has to satisfy the $\bf case\, 1$ or $\widetilde{u}\in
V_2(i_0-1,\e).$ Then, we define $Z_1(\widetilde{u})$ the
corresponding vector field, and we get
\begin{equation}\label{cR4}\bigl\langle\partial J(u), Z_1\bigr\rangle \leq -c\Bigl(\sum_{i<
i_0}\frac{1}{\lambda_i^{n-2}}+\frac{\bigl|\nabla
K(a_i)\bigr|}{\lambda_i}+\sum_{i \neq j,\,i,\,<
i_0}\e_{ij}\Bigr)+O\bigl(\sum_{i< i_0,\,j\geq
i_0}\frac{1}{\bigl(\lambda_i\lambda_j\bigr)^{\frac{n-2}{2}}}\bigr).\end{equation}
Observe that$$\lambda_i=o\bigl(\lambda_j\bigr),\,\forall\,\,i<
i_0,\,\forall\,\,j\geq i_0.$$Thus (\ref{cR4}) becomes
\begin{equation}\label{cR5}\bigl\langle\partial J(u), Z_1\bigr\rangle \leq -c\Bigl(\sum_{i<
i_0}\frac{1}{\lambda_i^{n-2}}+\frac{\bigl|\nabla
K(a_i)\bigr|}{\lambda_i}+\sum_{i \neq j,\,i,\,<
i_0}\e_{ij}\Bigr).\end{equation} Observe that all the
$\frac{1}{\lambda_i^{n-2}}'$s, $i_0\leq i\leq p,$ appear in the
upper bound (\ref{cR5}). Thus (\ref{cR5}) becomes
\begin{equation}\label{cR7}\begin{aligned}\bigl\langle\partial J(u), Z_1\bigr\rangle &\leq
-c\Bigl(\sum_{i=1}^p\frac{1}{\lambda_i^{n-2}} +\sum_{i<
i_0}\frac{\bigl|\nabla K(a_i)\bigr|}{\lambda_i}+\sum_{i \neq
j,\,i\,< i_0}\e_{ij}\Bigr)\\&\leq
-c\Bigl(\sum_{i=1}^p\frac{1}{\lambda_i^{n-2}} +\sum_{i<
i_0}\frac{\bigl|\nabla K(a_i)\bigr|}{\lambda_i}+\sum_{i \neq
j,\,1\leq i,\,j\leq p}\e_{ij}\Bigr).\end{aligned}\end{equation} Now,
arguing as in the proof of lemma \ref{l1}, we
get\begin{equation}\label{cR8}\bigl\langle\partial J(u),
Z_1+\sum_{i\geq i_0}\overline{X}_i\bigr\rangle \leq
-c\Bigl(\sum_{i=1}^p\frac{1}{\lambda_i^{2}} +\frac{\bigl|\nabla
K(a_i)\bigr|}{\lambda_i}+\sum_{i \neq j,\,1\leq i,\,j\leq
p}\e_{ij}\Bigr).\end{equation} The vector field
$\widetilde{W}_1^{1}:=Z_1+\sum_{i\geq i_0}\overline{X}_i$ satisfies
lemma \ref{l2}.\\ ${\bf Subcase\, 2.2}:$ $i_0= 1.$ Recall that
\begin{equation}\label{*1}\begin{aligned}
\bigl\langle \partial J(u), \sum_{i=1}^p\overline{X}_i\bigr\rangle
&\leq
-c\sum_{i=1}^p\Bigl[1-\Psi\bigl(\lambda_i|(a_i-y_{j_i})_{k_i}|\bigr)\Bigr]\frac{\bigl|\nabla
K(a_i)\bigr|}{\lambda_i}\\&+o\bigl(\sum_{i=1}^p\frac{\bigl|\nabla
K(a_i)\bigr|}{\lambda_i}+\frac{1}{\lambda_i^{n-2}}\bigr),\end{aligned}
\end{equation}where $k_i$ denotes the index such that $|(a_i-y_{j_i})_{k_i}|=\mathrm{max}_{1\leq k\leq
n}|(a_i-y_{j_i})_k|.$ Now, observe that in this ${\bf case\, 2}$, we
have\begin{equation}\label{*2}\frac{1}{\lambda_{i_1^{n-2}}}\leq
c\frac{1}{\delta^{\beta_j-1}}\frac{\bigl|\nabla
K(a_{i_1})\bigr|}{\lambda_{i_1}}\end{equation}
and\begin{equation}\label{*3}\Psi\bigl(\lambda_{i_1}|a_{i_1}-y_{j_{i_1}}|\bigr)=0.\end{equation}
Combining (\ref{*1}) and (\ref{*3}), we
obtain\begin{equation}\label{*4}\begin{aligned} \bigl\langle
\partial J(u), \sum_{i=1}^p\overline{X}_i\bigr\rangle &\leq
-c\Bigl(\frac{\bigl|\nabla K(a_{i_1})\bigr|}{\lambda_{i_1}}
+\sum_{i=1}^p\Bigl[1-\Psi\bigl(\lambda_i|(a_i-y_{j_i})_{k_i}|\bigr)\Bigr]\frac{\bigl|\nabla
K(a_i)\bigr|}{\lambda_i}\Bigr)\\&+o\bigl(\sum_{i=1}^p\frac{\bigl|\nabla
K(a_i)\bigr|}{\lambda_i}+\frac{1}{\lambda_i^{n-2}}\bigr).\end{aligned}
\end{equation}From the observation (\ref{*2}), we can make appear
$\frac{1}{\lambda_{i_1}^{n-2}}$ in the upper bound (\ref{*4}), and
then all the $\frac{1}{\lambda_{i}^{n-2}}'$s, $1\leq i\leq p.$ Then,
we obtain\begin{equation}\label{*5}\begin{aligned} \bigl\langle
\partial J(u), \sum_{i=1}^p\overline{X}_i\bigr\rangle &\leq
-c\Bigl(\frac{\bigl|\nabla
K(a_{i_1})\bigr|}{\lambda_{i_1}}+\sum_{i=1}^p\frac{1}{\lambda_{i}^{n-2}}
+\sum_{i=1}^p\Bigl[1-\Psi\bigl(\lambda_i|(a_i-y_{j_i})_{k_i}|\bigr)\Bigr]\frac{\bigl|\nabla
K(a_i)\bigr|}{\lambda_i}\Bigr)\\&\quad+o\bigl(\sum_{i=1}^p\frac{\bigl|\nabla
K(a_i)\bigr|}{\lambda_i}+\frac{1}{\lambda_i^{n-2}}\bigr)\\&\leq
-c\Bigl(\frac{\bigl|\nabla
K(a_{i_1})\bigr|}{\lambda_{i_1}}+\sum_{i=1}^p\frac{1}{\lambda_{i}^{n-2}}+\sum_{i\neq
j,\,1\leq i,\,j\leq p}\e_{ij}
\\&\quad\quad\quad+\sum_{i=1}^p\Bigl[1-\Psi\bigl(\lambda_i|(a_i-y_{j_i})_{k_i}|\bigr)\Bigr]\frac{\bigl|\nabla
K(a_i)\bigr|}{\lambda_i}\Bigr)\\&\quad+o\bigl(\sum_{i=1}^p\frac{\bigl|\nabla
K(a_i)\bigr|}{\lambda_i}+\frac{1}{\lambda_i^{n-2}}\bigr).\end{aligned}
\end{equation}Arguing as in the proof of the estimate (\ref{*10}),
the estimate (\ref{*5}) becomes
\begin{equation}\label{*6}\begin{aligned} \bigl\langle
\partial J(u), \sum_{i=1}^p\overline{X}_i\bigr\rangle &\leq
-c\Bigl(\sum_{i=1}^p\frac{1}{\lambda_{i}^{n-2}}+\frac{\bigl|\nabla
K(a_{i})\bigr|}{\lambda_{i}}+\sum_{i\neq j,\,1\leq i,\,j\leq
p}\e_{ij} \Bigr).\end{aligned}
\end{equation}The vector field
$\sum_{i=1}^p\overline{X}_i$ satisfies lemma \ref{l2}\\The
pseudo-gradient $\widetilde{W}_1$ required in lemma \ref{l2} will be
defined by convex combination of $Y_1,$ $\widetilde{W}_1^1$ and
$\sum_{i=1}^p\overline{X}_i.$\\

Observe that the variation of the maximum of the $\lambda_i'$s
occurs only under the condition
$$\lambda_i^{\alpha_{j_i}}|a_i-y_{j_i}| \leq \delta,\,\forall\,1\leq
i\leq p,$$ for $\delta> 0$ a fixed constant small enough. In this
case all the $\lambda_i'$s increase and go to $+\infty$ along the
flow-lines generated by $\widetilde{W}_1.$
\begin{lemma}\label{l3}\,In
$V_3(p,\e,d_0)$, there exists a pseudo-gradient $\widetilde{W}_3$ so
that the following holds: There is a constant $c > 0$ independent of
$u \in V_3(p,\e,d_0)$ so that$$\bigl\langle
\partial J(u), \widetilde{W}_3(u)\bigr\rangle\leq
-c\Bigl(\sum_{i=1}^p\frac{1}{\lambda_i^{n-2}}+\frac{\bigl|\nabla
K(a_i)\bigr|}{\lambda_i}+ \sum_{i \neq j}\e_{ij}\Bigr).$$\end{lemma}
\pf For each critical point $y_k$ of $K$, we set $B_k:=\{j,\,\,a_j
\in B(y_k, \eta)\}.$ Without loss of generality, we can assume
$y_1,\dots,\,y_q$ are the critical points such that
$\mathrm{card}(B_k)\geq 2,\, \forall \,k=1,\dots, q.$ Let $\chi$ be
a smooth cut-off function such that $\chi \geq0,\,\chi
=0\,\mbox{if}\,t\leq \gamma,\,\mbox{and}\,\chi=1\,\mbox{if}\,t\geq
1,$ where $\gamma$ is a small constant. Set
$\overline{\chi}(\lambda_j)=\sum_{i\neq j,\\i,j\in B_k}\,\chi
(\lambda_j/\lambda_i).$ Define$$W_3= -\sum_{k=1}^q\,\sum_{j\in
B_k}\,\alpha_j\overline{\chi}(\lambda_j)\lambda_j\frac{\partial P
\delta_j}{\partial\lambda_j}.$$Using proposition \ref{p1}, we derive
that\begin{equation*}\label{}
\begin{aligned}
\bigl\langle \partial J(u), W_3\bigr\rangle &=
2c_2J(u)\sum_{k=1}^q\sum_{j\in
B_k}\alpha_j\overline{\chi}(\lambda_j)\Biggl[\frac{n-2}{2}\alpha_j\frac{H(a_j,a_j)}{\lambda_j^{n-2}}\\&
\qquad\qquad\qquad\qquad\qquad\qquad +\sum_{i\neq
j}\alpha_i(\lambda_j\frac{\partial
\e_{ij}}{\partial\lambda_j}+\frac{n-2}{2} \frac{H(a_i,a_j)
}{(\lambda_i\lambda_j)^{\frac{n-2}{2}}})\Biggr]
\\&+ o\Bigl(\sum_{i \neq k}\e_{ik} +
\sum_{i=1}^p\frac{1}{\lambda_i^{n-2}}\Bigr)+\sum_{k=1}^q\,\sum_{i\in
B_k,\,\overline{\chi}(\lambda_i)\neq
0}\bigl(\mathrm{if\,\,}\lambda_i|a_i-y_{j_i}|\geq
C\bigr)o\bigl(\frac{\bigl|\nabla K(a_i)\bigr|}{\lambda_i}\bigr).
\end{aligned}
\end{equation*}
For $j\in B_k,$ with $k\leq q,$ if $\overline{\chi}(\lambda_j)\neq
0,$ then there exists $i\in B_k$ such that
$\lambda_j^{-1}=o(\lambda_i^{-1})$ or
$\lambda_j^{2-n}=o(\e_{ij})\,(\,\mbox{for}\,\,\eta\,\,\mbox{small\,
enough}).$ Furthermore, for $j\in B_k,$ if $i\not\in B_k$ or $i\in
B_k,$ with $\gamma < \lambda_i/\lambda_j< 1/\gamma,$ then we have
$\lambda_j\frac{\partial
\e_{ij}}{\partial\lambda_j}=-\e_{ij}(1+o(1)).$ In the case where
$i\in B_k$ and $\lambda_i/\lambda_j\leq \gamma,$ we have
$\overline{\chi}(\lambda_j)-\overline{\chi}(\lambda_i)\geq 1$ and
$\lambda_i\frac{\partial \e_{ij}}{\partial\lambda_i}\leq
-\lambda_j\frac{\partial \e_{ij}}{\partial\lambda_j}.$ Thus
$$\overline{\chi}(\lambda_j)\lambda_j\frac{\partial
\e_{ij}}{\partial\lambda_j}+\overline{\chi}(\lambda_i)\lambda_i\frac{\partial
\e_{ij}}{\partial\lambda_i}\leq \lambda_j\frac{\partial
\e_{ij}}{\partial\lambda_j}=-\e_{ij}(1+o(1))$$Thus we derive
that\begin{equation*}\label{}
\begin{aligned}\bigl\langle \partial J(u), W_3\bigr\rangle &\leq -c\sum_{k=1}^q\,\sum_{j\in
B_k,\,\overline{\chi}(\lambda_j)\neq 0}(\sum_{i \neq j}\e_{ij})+
o\Bigl(\sum_{i \neq k}\e_{ik} +
\sum_{i=1}^p\frac{1}{\lambda_i^{n-2}}+\frac{\bigl|\nabla
K(a_i)\bigr|}{\lambda_i}\Bigr)\\&\leq -c\sum_{k=1}^q\,\sum_{j\in
B_k,\,\overline{\chi}(\lambda_j)\neq
0}(\frac{1}{\lambda_j^{n-2}}+\sum_{i \neq j}\e_{ij})+ o\Bigl(\sum_{i
\neq k}\e_{ik} +
\sum_{i=1}^p\frac{1}{\lambda_i^{n-2}}\Bigr)\\&\quad+\sum_{k=1}^q\,\sum_{i\in
B_k,\,\overline{\chi}(\lambda_i)\neq
0}\bigl(\mathrm{if\,\,}\lambda_i|a_i-y_{j_i}|\geq
C\bigr)o\bigl(\frac{\bigl|\nabla
K(a_i)\bigr|}{\lambda_i}\bigr).\end{aligned}
\end{equation*}
Observe that $\{j\in B_k,\,\overline{\chi}(\lambda_j)=0\}$ contains
at most one index. Thus we obtain\begin{equation}\label{c6}
\begin{aligned}\bigl\langle\partial J(u), W_3\bigr\rangle \leq&
-c\bigl(\sum_{k=1}^q\,\sum_{j\in
B_k,\,\overline{\chi}(\lambda_j)\neq
0}\frac{1}{\lambda_j^{n-2}}+\sum_{i \neq j,\,j\in
\cup_{k=1}^qB_k}\e_{ij}\bigr)+ o\Bigl(\sum\limits_{i \neq k}\e_{ik}
+
\sum_{i=1}^p\frac{1}{\lambda_i}\Bigr)\\&\,\,\,+\sum_{k=1}^q\,\sum_{i\in
B_k,\,\overline{\chi}(\lambda_i)\neq
0}\bigl(\mathrm{if\,\,}\lambda_i|a_i-y_{j_i}|\geq
C\bigr)o\bigl(\frac{\bigl|\nabla
K(a_i)\bigr|}{\lambda_i}\bigr).\end{aligned}
\end{equation} This upper
bound does not contains all the indices. We need to add some terms.
Let$$\lambda_{i_0}=\mbox{inf}\{\lambda_i,\,i=1,\dots,\,p\}.$$Two
cases may occur:\\{\bf Case 1:} If $i_0\in
\cup_{k=1}^qB_k,\,\mbox{with}\,\overline{\chi}(\lambda_{i_0})\neq
0$, then we can make appear in the last upper bound
$\frac{1}{\lambda_{i_0}^{n-2}},$ and therefore all the
$\frac{1}{\lambda_i^{n-2}},$ and so $\e_{ik},\,1\leq i,\,k\leq
p.$\\{\bf Case 2:} $i_0 \in \{i\in
\cup_{k=1}^qB_k,\,\overline{\chi}(\lambda_i)=0\}\cup(\cup_{k=1}^qB_k)^c$.
In this case, we define$$D=(\{i\in
\cup_{k=1}^qB_k,\,\overline{\chi}(\lambda_j)=0\}\cup(\cup_{k=1}^qB_k)^c)\cap
\{1\leq i\leq p,\,\lambda_i/\lambda_{i_0}< 1/\gamma\}.$$It is easy
to see that for $i,\,j\in D,\,i\neq j,$ we have $a_i\in
B(y_{k_i},\eta)\,\,\mbox{and}\,\,a_j\in
B(y_{k_j},\eta)\,\,\mbox{with}\,\,k_i\neq k_j$. Let$$u_1=\sum_{i\in
D}\alpha_iP\delta_{(a_i,\,\lambda_i)}.$$$u_1$ has to satisfy one of
the two above cases, that is, $u_1\in
V_i(\mathrm{card}(D),\,\e,\,d_0)\,\mbox{for}\,i=1,\,2.$ Thus we can
apply the associated vector field which we will denote $W^{'}_3,$
and we have the following estimate:$$\bigl\langle\partial J(u),
W^{'}_3\bigr\rangle \leq -c\bigl(\sum_{i\in
D}\frac{1}{\lambda_i^{n-2}}+\frac{\bigl|\nabla
K(a_i)\bigr|}{\lambda_i}+\sum_{i \neq j,\,i,j\in D}\e_{ij}\bigr)+
O\Bigl(\sum\limits_{r\not\in D,\,k\in D}\e_{rk} + \sum_{i\not\in
D}\frac{1}{\lambda_i^{n-2}}\Bigr).$$Observe that for $k\in D$ and
$r\not\in D,$ we have either $r\in
\cup_{k=1}^qB_k,\,\overline{\chi}(\lambda_r)\neq 0$ ( in this case
we have $\e_{kr}$ in the upper bound (\ref{c6}) ) or no, and in this
last case we observe that $a_i\in
B(y_{j_i},\eta),\,\,\mbox{for}\,i=r,\,k,\,\,\mbox{with}\,\,j_r\neq
j_k$. Thus$$\e_{kr}\leq
\frac{c}{(\lambda_k\lambda_r)^{\frac{n-2}{2}}}\leq
\frac{c\gamma}{\lambda_{i_0}^{n-2}}=o(\frac{1}{\lambda_{i_0}^{n-2}}).$$We
get the same observation for $\lambda_i^{2-n},\,i\not\in D.$ Now we
define
$$Y_2=CW_3+W^{'}_3,$$where $C$ is a large positive constant. We obtain$$\bigl\langle\partial
J(u), Y_2\bigr\rangle \leq
-c\bigl(\sum_{i=1}^p\frac{1}{\lambda_i^{n-2}}+\sum_{i \neq
j}\e_{ij}\bigr)+o\bigl( \sum_{i=1}^p\frac{\bigl|\nabla
K(a_i)\bigr|}{\lambda_i}\bigr).$$ We define $W^{''}_2$ as a convex
combination of $W_3$ and $Y_2.$ Then the pseudo-gradient$$\widetilde{W}_3:=W^{''}_2+\sum_{i=1}^p\,\overline{X}_i$$satisfies
the claim of lemma \ref{l3}.\begin{lemma}\label{l4}\, In
$V_4(p,\e,d_0)$, there exists a pseudo-gradient $\widetilde{W}_4$ so
that the following holds: There is a constant $c > 0$ independent of
$u \in V_4(p,\e,d_0)$ so that$$\bigl\langle
\partial J(u),\widetilde{W}_4(u)\bigr\rangle\leq
-c\Bigl(\sum_{i=1}^p \frac{1}{\lambda_i^{n-2}}+\frac{|\nabla
K(a_i)|}{\lambda_i}+ \sum_{i \neq j}\e_{ij}\Bigr).$$\end{lemma}\pf
Without loss of generality, we suppose $\lambda_1\leq \dots\leq
\lambda_p.$\\We denote by $i_1$ the index satisfying $a_{i_1}\not
\in \cup_{\nabla K(y)=0}B(y, \eta)$ and $a_i \in B(y_{j_i},
\eta),\,\forall\,i<i_1.$ Let
$$\widetilde{u} = \sum_{i<i_1}\alpha_iP\delta_i.$$ Observe that
$\widetilde{u}\in V_i(i_1-1, \e),\,i=1,\,2\,\,\text{or}\,\,3.$ Then
we define $Z_4{'}(\widetilde{u})$ the corresponding vector field and
we have
$$\bigl\langle\partial J(u), Z_4{'}(\widetilde{u})\bigr\rangle\leq -c\Bigl(\sum_{i<i_1}\frac{1}{\lambda_i^{n-2}}+\frac{|\nabla
K(a_i)|}{\lambda_i} + \sum_{i \neq
j,\,i,j<i_1}\e_{ij}\Bigr)+O\Bigl(\sum_{i<i_1,\,j\geq
i_1}\e_{ij}+\sum_{i\geq i_1}\frac{1}{\lambda_i^{n-2}}\Bigr).$$Let
now
$$Z_4 =\frac{1}{\lambda_{i_1}}\frac{\nabla K(a_{i_1})}{|\nabla
K(a_{i_1})|}\frac{\partial P\delta_{a_{i_1},
\lambda_{i_1}}}{\partial a_{i_1}}-M_3\sum_{i\geq
i_1}2^i\lambda_i\frac{\partial P\delta_{a_i, \lambda_i}}{\partial
\lambda_i},$$ where $M_3>0$ is a fixed constant large enough. From
propositions \ref{p1} and \ref{p2}, we obtain, since $\nabla
K(a_{i_1})\neq 0,$\begin{equation*}\label{}
\begin{aligned}\bigl\langle
\partial J(u), Z_4(u)\bigr\rangle&\leq
\frac{-c}{\lambda_{i_1}}+O\Bigl( \sum_{j\neq
i_1}\lambda_j|a_{i_1}-a_j|\e_{i_1j}^3\Bigr)-M_3c\sum_{i\geq
i_1,j\neq i}\e_{ij}.\end{aligned}
\end{equation*}Observe that $\lambda_j|a_{i_1}-a_j|\e_{i_1j}^3=O(\e_{i_1j}),\,\forall\,j\neq i_1.$ Thus$$\bigl\langle
\partial J(u), Z_4(u)\bigr\rangle\leq
\frac{-c}{\lambda_{i_1}}+O\Bigl( \sum_{j\neq
i_1}\e_{i_1j}\Bigr)-M_3c\sum_{i\geq i_1,j\neq i}\e_{ij}.$$We choose
$M_3$ large enough so that $O\Bigl(\sum_{j\neq i_1}\e_{i_1j}\Bigr)$
is absorbed by $M_3c\sum_{i\geq i_1,\, j\neq i}\e_{ij}.$ We deduce
\begin{equation}\label{c7}\bigl\langle
\partial J(u), Z_4(u)\bigr\rangle\leq
-c\Bigl(\frac{1}{\lambda_{i_1}}+\sum_{i\geq i_1,\,i\neq
j}\e_{ij}\Bigr).\end{equation}Also $\frac{1}{\lambda_{i_1}}$ makes
appear $\sum_{i\geq i_1}\frac{1}{\lambda_i}$ in the upper bound of
(\ref{c7}). Taking M a positive constant large enough and
let$$\widetilde{W}_4(u)=MZ_4+Z_4^{'}.$$Thus we derive$$\bigl\langle
\partial J(u), \widetilde{W}_4(u)\bigr\rangle\leq
-c\Bigl(\sum_{i=1}^p\frac{|\nabla
K(a_i)}{\lambda_i}+\frac{1}{\lambda_i^{n-2}}+\sum_{i \neq
j}\e_{ij}\Bigr).$$The claim of lemma \ref{l4} follows.\\\\
The pseudo-gradient $W_1,$ required in proposition \ref{p3}, will be
defined by a convex combination of the vector fields
$\widetilde{W}_1(u)$, $\widetilde{W}_2(u)$, $\widetilde{W}_3(u)$ and
$\widetilde{W}_4(u).$\vskip3mm\n \pfn {\bf proposition \ref{p4}.} We
will introduce some technical lemmas for the proof of proposition
\ref{p4}. Without loss of generality, we suppose $\lambda_1d_1 \leq
... \leq \lambda_pd_p.$ Let $c_1 > 0$ a fixed constant small enough.
We define $$I_2:= \{1\} \cup \{ 1 \leq i \leq p,\, \text{s.t}\,\,
c_1\lambda_kd_k \leq \lambda_{k-1}d_{k-1} \leq \lambda_kd_k, \,\,
\forall \, k \leq i\}.$$ In $I_2$, we order the $\lambda_i$'s :
$\lambda_{i_1} \leq \dots \leq \lambda_{i_s}$. For $c_2 > 0$ a fixed
constant small enough, we define$$I_{\lambda_{i_s}}:= \{i_s\} \cup
\{1 \leq k \leq s,\,\text{ s.t}\,\, c_2\lambda_{i_{j+1}} \leq
\lambda_{i_j} \leq \lambda_{i_{j+1}},\, \forall \, j \geq k \}.$$For
$u = \sum\limits_{i=1}^p\alpha_iP\delta_i \in V_b(p,\e)$, we
introduce the following condition: $\mathrm{for}\,i\in
\{1,\dots,p\},$
\begin{equation}\label{condition} \frac{1}{2^{p+1}}\sum_{k \neq i}\e_{ki} \leq
\sum_{j=1}^p\frac{H_{ij}}{(\lambda_i\lambda_j)^{\frac{n-2}{2}}}.\end{equation}
We divide the set $\{1, \dots,\, p\}$  into $T_1 \cup T_2$, where
$$T_1 = \{1 \leq i \leq p,\, \mathrm{ s.t} \,\, i \,\, \mathrm{satisfies} \,\, (\ref{condition})\,\}$$and
$$T_2 = \{1,\dots, p\} \backslash T_1.$$
\begin{lemma}\label{l5}\,There exists a vector field $X_1$ such that
$$\bigl\langle\partial J(u), X_1 \bigr\rangle \leq
-c\Bigl(\frac{1}{\lambda_{i_s}} +
\sum_{i=1}^p\,\frac{1}{(\lambda_id_i)^{n-1}} + \sum_{k \in T_2, 1
\leq j \leq p } \e_{kj}+\sum\limits_{\,j\neq k;\,j\in T_1,\,1\leq
k\leq p\,}\e_{jk}^{\frac{n-1}{n-2}}\Bigr)+
o\Bigl(\sum_{i=1}^p\frac{1}{\lambda_i}\Bigr).$$\end{lemma}\pf We
define the vector field $Y_1^2$ by
$$Y_1^2:=\frac{1}{\lambda_{i_s}}\sum_{i\in
I_{\lambda_{i_s}}}\frac{\partial P\delta_i}{\partial
a_i}(-\alpha_i\nu_i).$$ From proposition \ref{p2}, we obtain
\begin{align}\label{upper1}
\begin{split}\bigl\langle\partial J(u), Y_1^2 \bigr\rangle &\leq
-c\Bigl(\frac{1}{\lambda_{i_s}} +
\frac{1}{(\lambda_{i_s}d_{i_s})^{n-1}} \Bigr)  +
\frac{1}{\lambda_{i_s}} O\Bigl(\sum_{k, j \in I_{\lambda_{i_s}}}
\lambda_k\lambda_j|a_k-a_j||\nu_k-\nu_j|\e_{kj}^{\frac{n}{n-2}}\Bigr)\\&+
\frac{1}{\lambda_{i_s}} O\Bigl(\sum_{k\in I_{\lambda_{i_s}}, j
\not\in I_{\lambda_{i_s}} }
\lambda_k\lambda_j|a_k-a_j|\e_{kj}^{\frac{n}{n-2}}\Bigr)\\& +
o\Bigl(\sum_{k \neq j}\e_{kj}^{\frac{n-1}{n-2}}\Bigr) +
o\Bigl(\sum_{k=1}^p\,\frac{1}{(\lambda_kd_k)^{n-1}}\Bigr).\end{split}
\end{align} Since $i_s$ belongs to $I_2$ and the term
$(\lambda_{i_s}d_{i_s})^{1-n}$ appear in the upper bound
$(\ref{upper1}),$ we can make appear the term $(\lambda_1d_1)^{1-n}$
in this upper bound, and so all the $(\lambda_id_i)^{1-n}.$ Observe
also, for $j \in T_1 \,\text{and}\,k\neq j,$ we have
$$\e_{kj}^{\frac{n-1}{n-2}} \leq
c\sum_{i=1}^p(\lambda_id_i)^{1-n}.$$Thus, after appearing
$\sum_{i=1}^p\frac{1}{(\lambda_id_i)^{n-1}},$ we have existence of
$\sum\limits_{j\neq k;\,j\in T_1,\,1\leq k\leq
p}\e_{jk}^{\frac{n-1}{n-2}}$ in the same upper bound. Observe that,
for $k, j \in I_{\lambda_{i_s}}$, $|\nu_k-\nu_j| = O(|a_k-a_j|)$.
So, we get$$\frac{1}{\lambda_{i_s}}
\lambda_k\lambda_j|a_k-a_j||\nu_k-\nu_j|\e_{kj}^{\frac{n}{n-2}} =
\frac{1}{\lambda_{i_s}}O(\e_{kj})
 = o(\frac{1}{\lambda_{i_s}}).$$ We are left for the
estimate of $\frac{1}{\lambda_{i_s}}
\lambda_k\lambda_j|a_k-a_j|\e_{kj}^{\frac{n}{n-2}}$ with $j \not\in
I_{\lambda_{i_s}}$ and $ k\in I_{\lambda_{i_s}} .$\\
If $k\in T_2\,\,\mathrm{or\,\,}j \in T_2$, we
get$$\frac{1}{\lambda_{i_s}}\lambda_k\lambda_j|a_k-a_j|\e_{kj}^{\frac{n}{n-2}}=O(\lambda_j|a_k-a_j|\e_{kj}^{\frac{n}{n-2}})
=O\Bigl((\frac{\lambda_j}{\lambda_k})^{\frac{1}{2}}
\e_{kj}^{\frac{n-1}{n-2}}\Bigr) =O\Bigl(\e_{kj}\Bigr).$$ If
$k,\,j\in T_1.$  In this case, we observe
that$$\frac{1}{\lambda_{i_s}}
\lambda_k\lambda_j|a_k-a_j|\e_{kj}^{\frac{n}{n-2}} =
o\Biggl(\sum_{i=1}^p \frac{1}{(\lambda_id_i)^{n-1}}\Biggr).$$As a
conclusion of the last observations, we obtain
\begin{align*}\begin{split}\bigl\langle\partial J(u), Y_1^2 \bigr\rangle \leq&
-c\Bigl(\frac{1}{\lambda_{i_s}}
+\sum_{i=1}^p\,\frac{1}{(\lambda_id_i)^{n-1}}+\sum\limits_{\,j\neq
k;\,j\in T_1,\,1\leq k\leq p\,}\e_{jk}^{\frac{n-1}{n-2}}\Bigr)
\\&+ O\Bigl(\sum_{k \in T_2, 1 \leq j \leq p }
\e_{kj}\Bigr).\end{split}\end{align*} Let $T_2  = \{i_1, \dots,
i_r\}, \, \text{with} \,\, \lambda_{i_1} \leq \dots \leq
\lambda_{i_r}$. We define
$$Y_2^2:= -
\sum_{k=1}^r2^k\alpha_{i_k}\lambda_{i_k}\frac{\partial P
\delta_{i_k}}{
\partial\lambda_{i_k}}.$$
Using proposition \ref{p1}, we obtain
\begin{equation}\label{c8}\begin{aligned}\bigl\langle\partial J(u), Y_2^2
\bigr\rangle=&2c_2J(u)\sum_{k=1}^r\,\Bigl[\sum_{j\neq
i_k}2^k\alpha_j\alpha_{i_k}\lambda_{i_k}\frac{\partial \e_{ji_k}}{
\partial\lambda_{i_k}}\\&{}\qquad\quad\qquad\,\,+\frac{n-2}{2}\sum_{j=1}^p2^k\alpha_j\alpha_{i_k}\frac{H(a_j,
a_{i_k})}{(\lambda_j\lambda_{i_k})^{\frac{n-2}{2}}})\Bigr]\times
\bigl(1+o(1)\bigr)\\&{}+o\Bigl(\sum_{i \neq j,\,i\in T_2}\e_{ij} +
\sum_{k=1}^{p}\frac{1}{(\lambda_kd_k)^{n-1}}+\frac{1}{\lambda_k}\Bigr).
\end{aligned}\end{equation}
We have$$\lambda_{i}\frac{\partial \e_{ij}}{
\partial\lambda_{i}}=-\frac{n-2}{2}\e_{ij}\bigl(1-2\frac{\lambda_j}{\lambda_i}\e_{ij}^{\frac{2}{n-2}}\bigr),\,\forall\,\,i\neq
j.$$ Observe also
that$$\frac{\lambda_j}{\lambda_i}\e_{ij}^{\frac{2}{n-2}}=o(1),\,\forall\,\,i\neq
j,\,j\in T_1.$$ Thus we derive
\begin{equation}\label{c9}
\lambda_{i}\frac{\partial \e_{ij}}{
\partial\lambda_{i}}\leq -\frac{3(n-2)}{8}\e_{ij},\,\forall\,\,i\neq
j,\,j\in T_1.
\end{equation}
However, for $j\in T_2$ and $\lambda_j\leq \lambda_i,$ we obtain
\begin{equation}\label{c10}
2\lambda_{i}\frac{\partial \e_{ij}}{
\partial\lambda_{i}}+\lambda_{j}\frac{\partial \e_{ij}}{\partial\lambda_{j}}\leq -\frac{3(n-2)}{8}\e_{ij},\,\forall\,\,i\neq j.
\end{equation} Combining (\ref{c8}), (\ref{c9}) and (\ref{c10}), we
derive
that\begin{equation}\label{c11}\begin{aligned}\bigl\langle\partial
J(u), Y_2^2 \bigr\rangle &\leq (n-2)c_2J(u)\sum_{k\in
T_2}\Bigl[-\frac{3}{4}\sum_{j\neq
k}\e_{jk}+2^p\sum_{j=1}^p\frac{H(a_j,
a_{k})}{(\lambda_j\lambda_{k})^{\frac{n-2}{2}}})\Bigr]\\&{}+o\Bigl(
\sum_{k=1}^{p}\frac{1}{(\lambda_kd_k)^{n-1}}+\frac{1}{\lambda_k}\Bigr)\\&{}\leq
-\sum_{j\neq k,\,k\in T_2}\e_{jk}+o\Bigl(
\sum_{k=1}^{p}\frac{1}{(\lambda_kd_k)^{n-1}}+\frac{1}{\lambda_k}\Bigr).
\end{aligned}\end{equation}
Taking $m> 0$ a fixed constant large enough, the vector field
$$X_1 =
Y_1^2 + mY_2^2$$satisfies the claim of lemma \ref{l5}.
\begin{lemma}\label{l6}\, There exist a vector field $X_2$ such that
$$\bigl\langle\partial J(u), X_2\bigr\rangle \leq
-c\Bigl(\sum_{i=1}^p\frac{1}{\lambda_i}\Bigr) + O\Bigl(\sum_{j\in
T_2, 1 \leq k \leq p}\e_{kj}\Bigr) +
O\Bigl(\sum_{i=1}^p\frac{1}{(\lambda_id_i)^{n-1}}\Bigr).$$\end{lemma}\pf
Without loss of generality, we suppose $\lambda_1 \leq \dots \leq
\lambda_p$. We define, for $M> 0$ a fixed constant large enough,
$$I_1^{'}:= \{ 1 \leq i \leq p,\,  \text{ s.t}\,\, |a_i-a_1| \geq
\frac{2}{M}d_1\},$$
$$I_1^{''}:= \{ i \not\in I_1^{'},\,\text{s.t}\,\, \exists (i_1,\dots,\,
i_r),\,\text{with} \, i_1=i, \,i_r\in I_1^{'},\, \text{and}\,\,
|a_{i_{k-1}}-a_{i_{k}}| < \frac{d_1}{pM},\, \forall \,\, k \leq
r\}$$and$$I_1:= \{1, \dots,\, p\}\backslash \{I_1^{'} \cup
I_1^{''}\}.$$Observe that, for $k\in I_1$ and $j\not\in I_1,$ we
have$$|a_k-a_j| \geq \frac{1}{pM}d_1.$$Let us define, for $c_2> 0$ a
fixed constant small enough,$$I_{\lambda_1}:= \{1\} \cup \{1 \leq j
\leq p,\, \text{s.t}\,\, c_2\lambda_{i_k} \leq \lambda_{i_{k-1}}
\leq \lambda_{i_k},\, \forall \, k \leq j \}.$$ We set
$$X_2:=\frac{1}{\lambda_1}\sum_{i\in I_1 \cap
I_{\lambda_1}} \frac{\partial P\delta_i}{\partial
a_i}(-\alpha_i\nu_i).$$ Observe that $d_i \sim d_1, \, \text{for}
\,\,i \in I_1\cap I_{\lambda_1}.$ From proposition \ref{p2}, we
obtain\begin{align*}\label{}
\begin{split}
\bigl\langle\partial J(u), X_2 \bigr\rangle &\leq  -
c\Bigl(\frac{1}{\lambda_1}\Bigr)  + \frac{1}{\lambda_1}
O\Bigl(\sum_{k, j \in I_{\lambda_1} \cap I_1}
\lambda_k\lambda_j|a_k-a_j||\nu_k-\nu_j|\e_{kj}^{\frac{n}{n-2}}\Bigr)\\&
+ \frac{1}{\lambda_1} O\Bigl(\sum_{k\in I_{\lambda_1} \cap I_1, j
\not\in I_{\lambda_1} \cap I_1 }
\lambda_k\lambda_j|a_k-a_j|\e_{kj}^{\frac{n}{n-2}}\Bigr)\\&+
o\Bigl(\sum_{k \neq j}\e_{kj}^{\frac{n-1}{n-2}}\Bigr) +
o\Bigl(\sum_{k=1}^p\,\frac{1}{(\lambda_kd_k)^{n-1}}\Bigr).
\end{split}
\end{align*}
Observe that, for $k, j \in I_{\lambda_1} \cap I_1$, $|\nu_k-\nu_j|
= O(|a_k-a_j|).$ From this, we deduce
$$\frac{1}{\lambda_1}
\lambda_k\lambda_j|a_k-a_j||\nu_k-\nu_j|\e_{kj}^{\frac{n}{n-2}} =
\frac{1}{\lambda_1}O(\e_{kj})=o(\frac{1}{\lambda_1}).$$ Now, we need
to estimate the quantity $\frac{1}{\lambda_1}
\lambda_k\lambda_j|a_k-a_j|\e_{kj}^{\frac{n}{n-2}},$
for $j \not\in I_{\lambda_1} \cap {I_1}\,\text{and}\,\, k\in I_{\lambda_1} \cap I_1.$ We have two cases:\\
{\bf Case 1: } $j \not\in I_1.$ In this case, we have $|a_k-a_j|
\geq \frac{1}{pM}d_1.$ From another side, since $k\in I_1,$ we
observe that $d_k \sim d_1$. Thus we deduce
$$\frac{1}{\lambda_1}
\lambda_k\lambda_j|a_k-a_j|\e_{kj}^{\frac{n}{n-2}}=
O\Bigl(\frac{1}{\lambda_1d_1}\e_{kj}\Bigr)=
O\Bigl(\frac{1}{(\lambda_1d_1)^{n-1}}\Bigr)+
o(\e_{kj}^{\frac{n-1}{n-2}}).$${\bf Case 2:} $j \in I_1.$ In this
case, we observe that $d_j \sim d_k \sim d_1 \,\, \mathrm{and}\,\,
|a_k-a_j| \leq \frac{4}{pM}d_1$.\\
If $j \in T_1.$ We use the fact that $H(a_j,\,a_j) \leq
c\frac{1}{d_j^{n-2}},\,\, \text{and}\,\, H(a_i,\,a_j) \leq
c\frac{1}{(d_id_j)^{\frac{n-2}{2}}}$ ( see \cite{50} ), then we
deduce\begin{align*}\label{}
\begin{split}
\frac{1}{\lambda_1}
\lambda_k\lambda_j|a_k-a_j|\e_{kj}^{\frac{n}{n-2}}&\leq\frac{1}{c_2^p}\lambda_j|a_k-a_j|\e_{kj}^{\frac{n}{n-2}}
\\&\leq \frac{c}{c_2^p}\lambda_j|a_k-a_j|\sum_{i=1}^p
\frac{1}{(\lambda_id_i)^{\frac{n}{2}}(\lambda_jd_j)^{\frac{n}{2}}}\\&
\leq \frac{c}{c_2^p}\frac{1}{M}\sum_{i=1}^p
\frac{1}{(\lambda_id_i)^{\frac{n}{2}}(\lambda_jd_j)^{\frac{n-2}{2}}}.
\end{split}
\end{align*}We choose $\frac{1}{c_2^p}\frac{1}{M} = o(1)$, and therefore
$$\frac{1}{\lambda_1}
\lambda_k\lambda_j|a_k-a_j|\e_{kj}^{\frac{n}{n-2}} =
o\Bigl(\sum_{i=1}^p
 \frac{1}{(\lambda_id_i)^{n-1}}\Bigr).$$
If $j \in T_2$, we easily have $$\frac{1}{\lambda_1}
\lambda_k\lambda_j|a_k-a_j|\e_{kj}^{\frac{n}{n-2}}=O\Bigl(\e_{kj}\Bigr).$$
We conclude that\begin{align*}\label{}
\begin{split}
\bigl\langle\partial J(u), X_2 \bigr\rangle \leq& -
c\Bigl(\frac{1}{\lambda_1}\Bigr) + O\Bigl(\sum_{j\in T_2, 1 \leq k
\leq p}\e_{kj}\Bigr) +
O\Bigl(\frac{1}{(\lambda_1d_1)^{n-1}}\Bigr)\\& +
o\Bigl(\sum_{k=1}^p\,\frac{1}{(\lambda_kd_k)^{n-1}}\Bigr).
\end{split}
\end{align*}
Such vector field $X_2$ satisfies the upper bound of lemma \ref{l6}.
Thus, for $m_1> 0$ a fixed constant large enough, the pseudo-gradient
$$W_2(u):= X_2 + m_1X_1$$satisfies the claim of
proposition \ref{p4}.\begin{corollary}\label{cor1}{\it\,Let $n \geq
4$. Assume that $K$ satisfies the condition $\mathbf{(f)}_\beta.$
Under the assumptions $\mathbf{(A_1)}$ and $\mathbf{(A_2)}$, The
critical points at infinity of $J$ in $V(p,\,\e),\,p\geq 1,$
correspond
to$$\sum_{j=1}^p\frac{1}{K(y_{i_j})^{\frac{n-2}{2}}}P\delta_{(y_{i_j},\,\infty)},$$
where $(y_{i_1},\,\dots,\,y_{i_p})\in \mathcal{C}_{\infty}.$
Moreover, such a critical point at infinity has an index equal to
$p-1+\sum_{j=1}^p\,n-\widetilde{i}(y_{i_j}).$}\end{corollary}\pf
Using theorem \ref{t2}, the only region where the $\lambda_i$'s are
unbounded is the one where each $a_i$ is close to a critical point
$y_{j_i}$ where $y_{j_i}\neq y_{j_k},$ for $i\neq k,$ and
$(y_{j_1},\dots,\,y_{j_p})\in \mathcal{C}_\infty.$ In this region,
arguing as in appendix 2 of \cite{6}, we can find a change of
variables$$(\,a_1,\,\dots,\,a_p,\,\lambda_1,\,\dots,\,\lambda_p\,)\mapsto
(\,\widetilde{a}_1,\,\dots,\,\widetilde{a}_p,\,\widetilde{\lambda}_1,\,\dots,\,\widetilde{\lambda}_p\,)
=:(\,\widetilde{a},\,\widetilde{\lambda}\,)
$$such that$$J(\sum_{i=1}^p\alpha_iP\delta_i +
\bar{v})=\frac{S_n^{\frac{2}{n}}\,\sum_{i=1}^p\alpha_i^2}{\bigl(
\sum_{i=1}^p\alpha_i^{\frac{2n}{n-2}}K(\widetilde{a}_i)\bigr)^{\frac{n-2}{n}}}
\bigl\{1 +c\cdot ^T\Lambda M \Lambda
\bigr\}=:\Psi\bigl(\alpha_1,\,\widetilde{a},\,\widetilde{\lambda}\bigr),$$
where $M:=M(y_{j_1},\dots,\,y_{j_p})$ is the matrix defined by
(\ref{matrice}),
$^T\Lambda:=(\frac{1}{(\lambda_1)^{\frac{n-2}{2}}},\dots,\,\frac{1}{(\lambda_p)^{\frac{n-2}{2}}}),$
$S_n:=\int_{\mathbb{R}^n}\delta_{o,\,1}^{\frac{2n}{n-2}}(x)dx$ and
$c$ is a positive constant. Observe that the function $\Psi$ admits
for the variables $\alpha_i'$s an absolute degenerate maximum with
one dimensional nullity space. Then the index of such critical point
at infinity is equal to
$p-1+\sum_{i=1}^p\,n-\widetilde{i}(y_{j_i}).$ The result of
corollary \ref{cor1} follows.
\section{Proof of the main result}\label{s4}
\def\theequation{4.\arabic{equation}}\makeatother
\setcounter{equation}{0}  \pfn {\bf theorem \ref{t1}}\,For technical
reasons, we introduce, for $\e_0> 0$ small enough, the following
neighborhood of $\Sigma^+:$
$$V_{\e_0}\bigl(\Sigma^+\bigr):=\{\,u\in \Sigma\,/\,\|u^-\|_{L^{\frac{2n}{n-2}}}< \e_0\,\},$$
where $u^-:=\mathrm{max}\bigl(0,\,-u\bigr).$ Recall that in theorem
\ref{t2} we construct a vector field $W$ defined in
$V(p,\,\e),\,\,p\geq 1.$ Outside $\bigcup\limits_{p\geq
1}V(p,\,\frac{\e}{2})$ we will use $-\partial J,$ and our global
vector field $Z$ will be built using a convex combination of $W$ and
$-\partial J.$ Arguing as in the proof of lemma 4.1 \cite{16}, we
can prove that $V_{\e_0}\bigl(\Sigma^+\bigr)$ is invariant under the
flow lines generated by $-\partial J,$ and therefore
$V_{\e_0}\bigl(\Sigma^+\bigr)$ is invariant under the flow lines
generated by $Z.$ We will prove the existence result by
contradiction. we suppose that $J$ has no critical points in
$V_{\e_0}\bigl(\Sigma^+\bigr).$ It follows from corollary \ref{cor1}
that the only critical points at infinity of $J$ in
$V_{\e_0}\bigl(\Sigma^+\bigr)$ correspond
to$$\sum_{j=1}^p\frac{1}{K(y_{i_j})^{\frac{n-2}{2}}}P\delta_{(y_{i_j},\,\infty)},\,\,p\geq
1, \,\mathrm{where}\,\, (y_{i_1},\,\dots,\,y_{i_p})\in
\mathcal{C}_{\infty}.$$ Such a critical point at infinity has an
index equal to $p-1+\sum_{j=1}^pn-\widetilde{i}(y_{i_j}).$\\

Given the pseudo-gradient $Z$ for $J$ on
$V_{\e_0}\bigl(\Sigma^+\bigr),$ we derive from the retraction
theorem 8.2 \cite{BR} that $V_{\e_0}\bigl(\Sigma^+\bigr)$ retracts
by deformation onto$$X_\infty:=\bigcup_{\tau_p\in
\mathcal{C}_\infty}\, W_u\bigl(\,(\tau_p)_\infty \,\bigr),$$where
$W_u\bigl(\,(\tau_p)_\infty\,\bigr )$ is the unstable manifold at
infinity associated to the critical point at infinity
$(\tau_p)_\infty.$ Since $V_{\e_0}\bigl(\Sigma^+\bigr)$ is
contractible, then we obtain
\begin{equation}\label{chi2}\chi\bigl(X_\infty\bigr)=\chi\Bigl(V_{\e_0}\bigl(\Sigma^+\bigr)\Bigr)=1.\end{equation}
Now, we call back the following fact: let $M$ be a finite cw-complex
in dimension $k,$ then the Euler-Poincar\'e characteristic of $M,$
$\chi\bigl(M\bigr),$ is given by
\begin{equation}\label{chi}\chi\bigl(M\bigr)=\sum_{i=1}^k\,(-1)^jn(j),\end{equation}
where $n(j)$ is the number of cells of dimension $j$ in $M$ ( see
\cite{Hatcher} ). We apply this fact to our situation, where the
cells of dimension an integer $j$ in $X_\infty$ are given by
$W_u\bigl(\,(\tau_p)_\infty\, \bigr)$ such that
$i\bigl(\,(\tau_p)_\infty\,\bigr)=j.$ According to (\ref{chi}), we
obtain\begin{equation}\label{chi1}\chi\bigl(X_\infty\bigr)=\sum_{\tau_p\in
\mathcal{C}_\infty}\,(-1)^{i\bigl(\tau_p\bigr)}.\end{equation}
Combining (\ref{chi1}) and (\ref{chi2}), we get$$\sum_{\tau_p\in
\mathcal{C}_\infty}\,(-1)^{i\bigl(\tau_p\bigr)}=1$$which contradicts
the assumption of our theorem. Thus there exists a critical point of
$J$ in $V_{\e_0}\bigl(\Sigma^+\bigr).$ Now, since $\e_0$ is small
enough, we derive by a standard argument that $u^-=0,$ and therefore
$u> 0\,\,\mathrm{in\,\,}\Omega.$ This finishes the proof of our
result.

\end{document}